\documentclass[12pt]{amsart}
\usepackage{multirow}
\usepackage{fullpage}
\usepackage{adjustbox}
\usepackage{float}
\usepackage{geometry}
\usepackage{amssymb,amsmath,amsthm}
\usepackage{bm}
\usepackage{bbm}
\usepackage{color}
\usepackage{enumerate}
\usepackage{fullpage}
\usepackage{mathtools}
\usepackage{makecell}

\newcommand{\bbC}{\mathbb{C}}

\newcommand{\bbF}{\mathbb{F}}
\newcommand{\bbZ}{\mathbb{Z}}

\newcommand{\calH}{\mathcal{H}}

\newcommand{\calW}{\mathcal{W}}
\newcommand{\calX}{\mathcal{X}}

\newcommand{\tr}{\mathrm{tr}}
\newcommand{\End}{\mathrm{End}}

\newcommand{\Hom}{\mathrm{Hom}}
\newcommand{\Ind}{\mathrm{Ind}}

\newcommand{\diag}{\mathrm{diag}}

\newcommand{\GL}{\mathrm{GL}}
\newcommand{\SL}{\mathrm{SL}}

\newcommand{\Ad}{\mathrm{Ad}}

\renewenvironment{proof}{\noindent{\scshape Proof.}}{\qed}

\theoremstyle{plain}
\newtheorem{theorem}{Theorem}[section]
\newtheorem{lemma}[theorem]{Lemma}
\newtheorem{proposition}[theorem]{Proposition}
\newtheorem{corollary}[theorem]{Corollary}

\newtheorem{definition}[theorem]{Definition}

\theoremstyle{definition}
\newtheorem{example}[theorem]{Example}
\newtheorem{remark}[theorem]{Remark}

\numberwithin{equation}{section}

\title[Multiplicity free induction]{Multiplicity free induction for the pairs \\
$(\GL_2\times\GL_{2},\mathrm{diag}(\GL_{2}))$ and $(\SL_3,\GL_2)$ over finite fields}

\author{Elias Depuydt}
\author{Maarten van Pruijssen}


\begin{document}

\begin{abstract}
We classify the irredible representations of $\GL_{2}(q)$ for which the induction to the product group $\GL_{2}(q)\times\GL_{2}(q)$, under the diagonal embedding, decomposes multiplicity free. It turns out that only the irreducible representations of dimensions $1$ and $q-1$ have this property.
We show that for $\GL_{2}(q)$ embedded into $\SL_{3}(q)$ via $g\mapsto\diag(g,\det g^{-1})$ none of the irreducible representations of $\GL_{2}(q)$ induce multiplicity free.
In contrast, over the complex numbers, the holomorphic representation theory of these pairs is multiplicity free and the corresponding matrix coefficients are encoded by vector-valued Jacobi polynomials. We show that similar results cannot be expected in the context of finite fields for these examples.
\end{abstract}

\maketitle


\section{Introduction}
This paper is a result of an investigation of multiplicity free induction of complex irreducible representations in the context of the matrix groups $\GL_{2}(q)\times\GL_{2}(q)$ and $\SL_{3}(q)$ over the finite field $\bbF_{q}$ with $q=p^{\ell}$ elements, in which the group $\GL_{2}(q)$ is embedded diagonally into $\GL_{2}(q)\times\GL_{2}(q)$ and via $g\mapsto\diag(g,\det g^{-1})$ into $\SL_{3}(q)$. The first pair is a Gelfand pair, which means that the trivial representation of $\GL_{2}(q)$ induces multiplicity free to $\GL_{2}(q)\times\GL_{2}(q)$, i.e.~the induced module decomposes into irreducible submodules with multiplicity at most one. The irreducible representations of $\GL_{2}(q)$ are of dimension $1$, $q-1,q$ or $q+1$.
We show in Section \ref{Section:product} that the irreducible representations of $\GL_{2}(q)$ of dimensions $1$ and $q-1$ induce multiplicity free to $\GL_{2}(q)\times\GL_{2}(q)$.

A triple $(G,H,\pi)$ with $G$ a finite group, $H\le G$ a subgroup and $\pi:H\to\GL(V)$ an irreducible complex $H$-representation that induces multiplicity free to $G$, is called a Gelfand triple.
Gelfand triples and their harmonic analysis have been investigated recently by Ceccherini-Silberstein, Scarabotti and Tolli \cite{MR4162267}. However, only a small number of such examples are known in the literature. 
This paper contributes a large new class of Gelfand triples arising from natural matrix group embeddings over finite fields. These examples are closely related to a well-understood classification over the complex numbers of large families of (holomorphic) Gelfand triples \cite{MR4642865}. While all the corresponding cases over finite fields could, in principle, be analyzed, the powerful geometric methods available over $\bbC$—notably those based on spherical varieties—do not yet have analogues in the finite field setting.

A striking observation is that the pair $(G, H)$ in a Gelfand triple $(G, H, \pi)$ need not itself be a Gelfand pair \cite[Thm.A.2]{MR4162267}, see Example \ref{example:GTnoGP}. This opens the possibility that, even if the pair $(\SL_{3}(q), \GL_{2}(q))$ is not a Gelfand pair, certain irreducible representations of $\GL_{2}(q)$ might still induce multiplicity-free to $\SL_{3}(q)$. In Section~\ref{Section:projective}, we show that this is not the case, no irreducible representations of $\GL_{2}(q)$ induce multiplicity-free to $\SL_{3}(q)$.

Note that, in contrast, the pairs of matrix groups that we consider in this paper, when taken with complex coefficients, are actually strong Gelfand pairs in the sense that each holomorphic irreducible representation of $\GL_{2}(\bbC)$ induces multiplicity free to both $\GL_{2}(\bbC)\times\GL_{2}(\bbC)$ and $\SL_{3}(\bbC)$. The harmonic analysis in the context over $\bbC$ was the main motivation for this paper, more background is given in Section \ref{Section:Background}.
In Section \ref{Section:SF} we discuss the the relevant parts of the harmonic analysis for the finite Gelfand triples $(G,H,\pi)$ to indicate that the relevant matrix coefficients cannot be described by families of vector-valued orthogonal polynomials for our examples, in contrast to the complex case.

\subsection{Notation} Throughout this paper $q=p^{\ell}$ is a prime power and $\GL_{2}(q)$ denotes the general linear group of $2\times 2$ matrices with coefficients in $\bbF_{q}$, the field with $q$ elements. Likewise $\SL_{3}(q)$ denotes the special linear group of $3\times 3$-matrices with coefficients in $\bbF_{q}$. We denote $\bbZ_{n}=\bbZ/(n)$ where $(n)=n\bbZ$.
The action of $G$ on itself by conjugation is denoted by $\Ad(g)(x)=gxg^{-1}$. The set of conjugacy classes is denoted by $G/\Ad(G)$. If the groups $G,H$ are fixed, we sometimes denote $\Ind V$ for the induced module $\Ind_{H}^{G}(V)$.

\section{Background and Motivation}\label{Section:Background}
Let $G$ be a finite group and $H\le G$ a subgroup.
In the case that all irreducible representations of $H$ induce multiplicity free to $G$, the convolution algebra of $H$-class functions on $G$ is commutative, and vice versa \cite[Thm.2, Cor.1 of Thm.7]{MR0340399} and $H$ is called a multiplicity free subgroup of $G$. For non-trivial examples we refer to \cite{MR4748964}, where the multiplicity free subgroups of $\SL_{2}(q)$ have been classified.
If there is an irreducible representation of $H$ that induces multiplicity free to $G$, then we say that $H$ admits a Gelfand triple. In general, not every subgroup admits a Gelfand triple.  
For example, in \cite{MR4815523} the subgroups of the symmetric groups $S_{n}$ that admit a Gelfand triple have been classified, for $n\ge66$.

The multiplicity free induction of an irreducible $H$-representation $\pi$ to $G$ can be detected by the commutativity of a suitable convolution subalgebra $I_{\pi}(G)$ of the $H$-class functions on $G$. If the algebra $I_{\pi}(G)$ is commutative, there is a rich harmonic analysis available including a spherical Fourier transform and a corresponding inversion formula \cite{MR1450744}.
We show that the algebra $I_{\pi}(G)$ is anti-isomorphic to the Hecke algebra $\widetilde{\calH}(G,H,\pi)$ of the triple $(G,H,\pi)$, see \cite[Def.2.1]{MR4162267}. The fact that the Hecke algebra encodes the multiplicity freeness can also be seen directly, as it is isomorphic to $\End_{G}(\Ind^{G}_{H}(\pi))$.

Spherical functions for the pair $(G,H)$ of type $\pi$, an irreducible representation of $H$, are matrix-valued functions that encode representations of $I_{\pi}(G)$. Spherical functions of type $\pi_{1}$, the trivial $H$-representation, are called zonal spherical functions.
The theory of spherical functions has a rich history, shaped by the contributions of many mathematicians. In the setting of locally compact groups, Godement’s seminal paper \cite{MR0052444} presents key ideas and cites earlier foundational work by Gelfand, Harish-Chandra, and others.
Recently, Blanco Villacorta, Pacharoni and Tirao, have given an exposition of this theory in the context of finite groups \cite{MR4794599}. 
Our discussion in Section~\ref{Section:SF} builds on the work of Gangolli and Varadarajan \cite{MR0954385}, which applies to general locally compact groups, and seeks to bridge these general ideas with those in \cite{MR4162267} by Ceccherini-Silberstein, Scarabotti, and Tolli, where a somewhat different convolution algebra is employed.
As a result, we observe that, for example Camporesi's results \cite{MR1450744}, apply to the harmonic analysis of the Gelfand triples in general and to our examples in particular.
The definition of spherical functions that we adopt in Section \ref{Section:SF} is more convenient when considering restrictions and inductions of representations. 

For symmetric pairs of compact Lie groups, the zonal spherical functions,  which can be viewed as members of a particular basis of the $H$-biinvariant functions on $G$, can be interpreted as Heckman-Opdam polynomials \cite[Ch.5]{MR1313912}. When considering a compact symmetric pair $(G,H)$, where an irreducible representation $\pi$ of $H$ induces multiplicity-free, the corresponding spherical functions — under mild assumptions — can be described using vector-valued orthogonal polynomials that share many properties with Heckman-Opdam polynomials. The reason is that the vector space spanned by the zonal spherical functions is a polynomial algebra and the vector space spanned by the the spherical functions of type $\pi$ is free and finitely generated over this algebra.
In such cases, the compact groups $G$ and $H$ can be complexified, and these results can be proved with the combinatorics of spherical varieties \cite[Thm.8.12]{MR4642865}.  
These techniques are not applicable to finite groups. However, there are still classes of spherical functions that are connected to orthogonal polynomials and special functions.

\begin{example}
    Consider Gelfand pairs $(G,P)$, where $G$ is a Chevalley group and $P$ a maximal parabolic subgroup. The zonal spherical functions for such pairs can be understood as Hahn- and Krawtchouk-polynomials \cite{MR0774056}. 
    We are not aware of any such description of the zonal spherical functions in the context of $\GL_{2}(q)$ diagonally embedded into $\GL_{2}(q)\times\GL_{2}(q)$. 
\end{example}

\begin{example}\label{example:GTnoGP}
    Let $G=\GL_{2}(q)$ and $U$ the subgroup of upper-triangular matrices with ones on the diagonal. The group $U$ is isomorphic to $\bbF_{q}$ and each non-trivial character $\psi:U\to\bbC^{\times}$ induces multiplicity free to $\GL_{2}(q)$, where the induced module decomposes into the sum of all irreducible $G$-modules of dimension $>1$, all with multiplicity one \cite[Thm.16.1]{P-S1983}. This makes $(G,U,\psi)$ a Gefland triple. However, $(G,U)$ is not a Gelfand pair.
    The induction of the trivial character of $U$ to $G$ is not multiplicity free, this follows for example from \cite[Thm.8.13]{P-S1983}, which implies that the irreducible $G$-modules of dimension $q+1$ occur with multiplicity two in the $G$-module induced by the trivial representation of $U$.
    The spherical functions of type $\psi$, a non-trivial character of $U$, are called Bessel functions, since they have properties similar to the classical Bessel functions, see e.g.~\cite{P-S1983}. 
    The vector space spanned by the spherical functions of type $\psi$ has dimension $q(q-1)$, the number of all irreducible representations of $G$ of dimension $>1$. 
    This space cannot be a free module over the algebra spanned by the zonal spherical functions, since the dimension of that algebra as a vector space is $(q-1)(q+1)^2$, the number of elements in $G/U$. 
\end{example}

These examples indicate a connection of classes of matrix coefficients on finite groups of Lie type with special functions and orthogonal polynomials. If we take the groups in these examples with coefficients in an algebraically closed field $K$, then we obtain spherical pairs: a Borel subgroup of $G(K)$ has an open orbit in the quotient variety $G(K)/H(K)$, see for example \cite{Timashev} for generalities and \cite{MR4642865} for relations with the multiplicity free representation theory. 

The two examples in this paper are among the smallest spherical pairs over~$\bbC$ with reductive subgroups. The classification of such pairs over~$\bbC$ is due to Krämer, Brion, and Mikityuk, see e.g.~\cite[\S10.2]{Timashev}, and has been extended to general algebraically closed fields by Knop and Röhrle \cite{KnopRohrle}. These classifications suggest a possible source for finding interesting Gelfand triples, but their realization over finite fields has not yet been systematically explored.

Our contribution is to show that $(\GL_{2}(q)\times \GL_{2}(q),\diag(\GL_{2}(q)),\pi)$
with $\pi$ irreducible of dimension $q-1$, is a Gelfand triple. This offers a new perspective: by starting from known spherical pairs over~$\bbC$, one can investigate which corresponding triples over finite fields define new multiplicity-free situations. This viewpoint opens a possible direction for developing new theory and examples.

Related results on tensor product decompositions of irreducible representations of $\GL_{2}(q)$ have recently been obtained by Gupta and Hassain \cite{GuptaHassain}, based on earlier work by Aburto, Johnson, and Pantoja \cite{MR2319497}. Their results appeared shortly before ours, but are based on different methods and arise in a different context and motivation.
Our results rely on elementary character theory, which also applies in the case of $\SL_{3}(q)$ using the character table in \cite{SimpsonFrameSutherland}. In fact, we show that the pair $(\SL_{3}(q),\GL_{2}(q))$ does not give rise to any Gelfand triple.


\section{Spherical functions}\label{Section:SF}

Let $G$ be a finite group and $\widehat{G}=\{\rho_{1},\ldots,\rho_{m}\}$, a complete set of pairwise inequivalent irreducible complex linear representations $\rho_{i}:G\to\GL(V_{i})$. 
For $\rho:G\to\GL(V)$ a linear representation, the map
$$P(G,\rho,\rho_{i})=\frac{\dim V_{i}}{|G|}\sum_{y\in G}\tr\rho_{i}(y^{-1})\rho(y)$$
is the projection of $V$ onto the $\rho_{i}$-isotypical subspace $V(\rho_{i})$ of $V$, i.e.~the sum of all the $G$-submodules of $V$ that are $G$-isomorphic to $V_{i}$. 

Let $H$ be a subgroup of $G$ with $\widehat{H}=\{\pi_{1},\ldots,\pi_{n}\}$ a complete set of pairwise inequivalent irreducible representations $\pi_{j}:H\to\GL(W_{j})$. Let $\pi_{1}$ be the trivial representation.
\begin{definition}\label{Def:SF}
The spherical function of type $\pi_{j}$ associated $\rho_{i}$ is the $\End(V_{i}(\pi_{j}))$-valued function on $G$ defined by
$$\Phi^{i}_{j}(g)=P(H,\rho_{i},\pi_{j})\rho_{i}(g)P(H,\rho_{i},\pi_{j}).$$
Spherical functions of type $\pi_{1}$, the trivial $H$-representation, are called zonal spherical function.
\end{definition}
The spherical functions of type $\pi_{j}$ do not take their values all in the same space due to various multiplicities $\dim\Hom_{H}(W_{i},V_{j})$.
Instead of looking at spaces of spherical functions, we use the spherical functions to study function spaces on $G$. 
Denote $\xi_{j}(h)=\dim W_{j}\tr\pi_{j}(h^{-1})$ which we view as a function on $G$ by declaring it zero outside of $H$.
Let $L(G)$ denote the space of $\bbC$-valued functions on $G$, which is an algebra for the convolution product
$$(f_{1}*f_{2})(x)=\frac{1}{|G|}\sum_{y\in G}f_{1}(y)f_{2}(y^{-1}x).$$
\begin{definition}
Let $L_{\pi_{j}}(G)$ be the space of functions $f\in L(G)$ for which $f=\xi_{j}*f=f*\xi_{j}$.
\end{definition}
The algebra $L_{\pi_{j}}(G)$ is a subalgebra of $L(G)$ and $f\mapsto \xi_{j}*f*\xi_{j}$ is a projection of $L(G)$ onto $L_{\pi_{j}}(G)$.
Denote $\langle f,\Phi^{i}_{j}\rangle=\frac{1}{|G|}\sum_{y\in G}f(y)\Phi^{i}_{j}(y)$.
The map $\gamma^{i}_{j}:L_{\pi_{j}}(G)\to\End(V(\pi_{j}))$ given by $\gamma^{i}_{j}(f)=\langle f,\Phi^{i}_{j}\rangle$ is a representation. Essentially every representation of $L_{\pi_{j}}(G)$ is of this form and this sets up a correspondence, up to equivalence, between the irreducible representations $\rho_{i}$ of $G$ with $\dim\Hom_{H}(W_{i},V_{j})>0$ and the non-trivial irreducible representations of $L_{\pi_{j}}(G)$ (cf.~\cite[\S1.3]{MR0954385}).

The analysis of $L_{\pi_{j}}(G)$ can be refined as follows. Let $I(G)=L(G)^{\Ad(H)}$, the subalgebra of $H$-class functions on $G$, i.e.~functions $f$ on $G$ for which
$$f(hgh^{-1})=f(g)\quad\mbox{for all $h\in H$ and $g\in G$.}$$
Let $I_{\pi_{j}}(G)=L_{\pi_{j}}(G)\cap I(G)$ which is equal to $I_{\pi_{j}}(G)=\xi_{\pi_{j}}*I(G)*\xi_{\pi_{j}}$. The algebra $L_{\pi_{j}}(G)$ can be recovered from $I_{\xi_{j}}(G)$ as  follows. Let $L_{\pi_{j}}(H)=\xi_{\pi_{j}}*L(H)*\xi_{\pi_{j}}$, the algebra of matrix coefficients on $H$ of the irreducible representation dual to $\pi_{j}$. As before we view functions on $H$ as functions on $G$ by extending them trivially outside $H$, so we can take convolution products of functions in $L(H)$ with functions in $L(G)$. The bilinear map $L_{\pi_{j}}(H)\times I_{\pi_{j}}(G)\to L_{\pi_{j}}(G)$ with $(a,f)\mapsto a*f$ extends to an isomorphism of algebras $L_{\pi_{j}}(H)\otimes I_{\pi_{j}}(G)\to L_{\pi_{j}}(G)$ \cite[L.1.3.10]{MR0954385}.

The non-trivial irreducible representations of $I_{\pi_{j}}(G)$ are also in a bijective correspondence, up to equivalence, with the irreducible representations $\rho_{j}$ of $G$ with $\dim\Hom_{H}(W_{i},V_{j})>0$. In view of the multiplicity freeness we record the following result \cite[Prop.1.3.15]{MR0954385}.

\begin{proposition}
The algebra $I_{\pi_{j}}(G)$ is commutative if and only if $I_{\pi_{j}}(G)$ is the center of $L_{\pi_{j}}(G)$. In this case we have $\dim\Hom_{H}(W_{j},V_{i})\le1$ for all $i=1,\ldots,m$.
Conversely, if $\dim\Hom_{H}(W_{j},V_{i})\le1$ for all $i=1,\ldots,m$, then $I_{\pi_{j}}(G)$ is commutative.
\end{proposition}
Note that $\dim\Hom_{H}(W_{j},V_{i})\le1$ for all $i=1,\ldots,m$ and $j=1,\ldots,n$ if and only if $I_{\pi_{j}}(G)$ is commutative for each $j=1,\ldots,n$.
In turn, since
$$L(G)^{\Ad(H)}=\bigoplus_{j=1}^{n}I_{\pi_{j}}(G),$$
this is equivalent to $L(G)^{\Ad(H)}$ being commutative. In this way we recover a result of Travis \cite[Cor.1 to Thm.7]{MR0340399}. A subgroup $H\le G$ with $L(G)^{\Ad(H)}$ commutative is called a multiplicity free subgroup.

\begin{definition}
A triple $(G,H,\pi_{j})$ is called a multiplicity free triple or a Gelfand triple if $I_{\pi_{j}}(G)$ is commutative.
\end{definition}

We can now calculate if $I_{\pi_{j}}(G)$ is commutative indirectly by calculating the multiplicities $\dim\Hom_{H}(W_{j},V_{i})$ using character theory. The space of class functions on $H$ is equipped with the inner product
\begin{equation}\label{FormulaInnerProduct}(\phi|\psi)_{H}=\frac{1}{|H|}\sum_{x\in H}\phi(x)\overline{\psi(x)}
\end{equation}
which is very useful to compute multiplicities. Indeed, if $\rho:H\to\GL(V)$ is a representation with character $\tr\rho$ then
\begin{equation}\label{FormulaMultiplicity}
(\tr\rho|\tr\pi_{j})_{H}=\dim\Hom_{H}(W_{j},V).
\end{equation}
For $\rho=\rho_{i}|_{H}$ the restriction to $H$ of the irreducible $G$-representation $\rho_{i}$ we denote
$$M^{i}_{j}=\dim\Hom_{H}(W_{j},V_{j}).$$
The matrix $(M^{i}_{j})_{i,j}$ is called the multiplicity table for $G$ and $H$. In the rows we see the irreducible constituents  and their multiplicities, of the restriction of $\rho_{i}$ to $H$. In the columns we see which irredicible $G$-representations contain copies of $\pi_{j}$, and how many, upon restriction to $H$. By Frobenius Reciprocity this tells us how $\Ind^{G}_{H}(W_{i})$ decomposes into irreducible $G$-modules. The $j$-th column has entries $\le 1$ if and only if $I_{\pi_{j}}(G)$ is commutative, i.e.~if and only if $(G,H,\pi_{j})$ is a multiplicity free triple. 

\begin{example}\label{Example:MFinduction}
Consider the permutation group $S_{4}$ with subgroup $\langle r\rangle$ with $r=(1~2~3)$, that is isomorphic to the cyclic group $C_{3}$. 
The group $S_{4}$ has five  conjugacy classes 
while $C_{3}$ has three. Their character values, taken from \cite[\S5.1,8]{Serre1977LinRep}, are recorded in Tables \ref{S_4chartable} and \ref{C_3chartable}.  The superscript indicates the number of elements in the corresponding conjugacy class.

\begin{table}[h]
    \centering
    \begin{minipage}{0.58\textwidth}
        \centering
        \begin{tabular}{|c |c c c c c|}
            \hline
            & Id$^{1}$ & $(1~2)^{6}$ & $(1~2~3)^{8}$ & $(1~2~3~4)^{6}$ & $(1~2)(3~4)^{3}$ \\
            \hline
            $\chi_1$ & $1$ & $1$ & $1$ & $1$ & $1$ \\
            $\chi_2$ & $1$ & $-1$ & $1$ & $-1$ & $1$ \\
            $\chi_3$ & $2$ & $0$ & $-1$ & $0$ & $2$ \\
            $\chi_4$ & $3$ & $1$ & $0$ & $-1$ & $-1$ \\
            $\chi_5$ & $3$ & $-1$ & $0$ & $1$ & $-1$ \\
            \hline
        \end{tabular}
        \caption{Character table of $S_4$.}
        \label{S_4chartable}
    \end{minipage}%
    \hfill
    \begin{minipage}{0.38\textwidth}
        \centering
        \begin{tabular}{|c |c c c|}
            \hline
            & Id$^{1}$ & $(1~2~3)^{1}$ & $(1~3~2)^{1}$ \\
            \hline
            $\psi_1$ & $1$ & $1$ & $1$ \\
            $\psi_2$ & $1$ & $e^{2\pi i/3}$ & $e^{4\pi i/3}$ \\
            $\psi_3$ & $1$ & $e^{4\pi i/3}$ & $e^{2\pi i/3}$ \\
            \hline
        \end{tabular}
        \caption{Character table of $C_3$}
        \label{C_3chartable}
    \end{minipage}
\end{table}

The restrictions of the characters of $S_{4}$ to $\langle r\rangle$ can be computed and we obtain the multiplicities that are recorded in Table \ref{S_4multtable}.
We observe that $\langle r\rangle\le S_{4}$ is a multiplicity free subgroup. 

\begin{table}[h]
\begin{center}
    \begin{tabular}{|c|c c c|}
    \hline
         & $\psi_1$ & $\psi_{2}$ & $\psi_{3}$ \\
         \hline
       $\mathrm{Res}\chi_1$  & $1$ & $0$ & $0$\\
       $\mathrm{Res}\chi_2$  & $1$ & $0$ & $0$\\
       $\mathrm{Res}\chi_3$  & $0$ & $1$ & $1$\\
       $\mathrm{Res}\chi_4$  & $1$ & $1$ & $1$\\
       $\mathrm{Res}\chi_5$  & $1$ & $1$ & $1$\\
       \hline
    \end{tabular}
    \caption{Multiplicity table of $(S_4,\langle~(1~2~3)~\rangle)$}\label{S_4multtable}
\end{center}
\end{table}
\end{example}

There is another convolution algebra related to the representation theory of the triple $(G,H,\pi_{j})$. Consider the space $L(G)\otimes\End(W_{j})$ of $\End(W_{j})$-valued functions.
Let $\widetilde{\mathcal{H}}(G,H,\pi_{j})$ denote the subspace of $L(G)\otimes\End(W_{j})$ consisting of functions $F:G\to\End(W_{j})$ satisfying the transformation behavior
$$F(h_{1}gh_{2})=\pi_{j}(h_{2})^{-1}F(g)\pi_{j}(h_{1})^{-1}$$
for all $h_{1},h_{2}\in H$ and $g\in G$.
The space $\widetilde{\mathcal{H}}(G,H,\pi_{j})$ is a convolution algebra for
$$(F_{1}*F_{2})(x)=\frac{1}{|G|}\sum_{y\in G}F_{1}(y^{-1}x)F_{2}(y),$$
cf.~\cite[\S2]{MR1450744}.
The algebra $\widetilde{\mathcal{H}}(G,H,\pi_{j})$ is called the Hecke algebra of the triple $(G,H,\pi_{j})$ \cite[Def.2.1]{MR4162267}.
Consider the the induced representation $\Ind^{G}_{H}\pi_{j}$ on the corresponding module $\Ind^{G}_{H}W_{j}$ of $W_{j}$-valued functions on $G$ with $\pi_{j}(h)f(g)=f(gh)$ for all $h\in H$ and $g\in G$, and the algebra of $G$-intertwiners $\End_{G}(\Ind^{G}_{H}W_{j})$. 
The map
$$\xi:\widetilde{\calH}(G,H,\pi_{j})\to\End_{G}(\Ind^{G}_{H}(W_{i}))$$
defined by $\xi(F)(f)(x)=\sum_{y\in G}	F(y^{-1}x)f(x)$ establishes an isomorphism of algebras \cite[Thm.2.1]{MR4162267}. We observe that $(G,H,\pi_{j})$ is a Gelfand triple if and only if the algebra $\End_{G}(\Ind^{G}_{H}(W_{j}))$ is commutative, where we use Schur's Lemma. In turn, this is equivalent to the Hecke algebra $\widetilde{\mathcal{H}}(G,H,\pi_{j})$ being commutative, where we use the isomorphism $\xi$. This is not surprising in view what follows.
From \cite[\S2]{MR1450744}, if $F\in\widetilde{\calH}(G,H,\pi_{j})$ then $f_{F}=\dim(W_{j})\tr F\in I_{\pi_{j}}(G)$. Conversely, if $f\in I_{\pi_{j}}(G)$ then
$$F_{f}(x)=\frac{1}{|H|}\sum_{y\in H}\pi_{j}(y)f(yx)$$
defines a function $F_{f}\in \widetilde{\calH}(G,H,\pi_{j})$. This gives an algebra anti-isomorphism $$\widetilde{\calH}(G,H,\pi_{j})\to I_{\pi_{j}}(G),\quad F\mapsto f_{F}$$
with inverse $f\mapsto F_{f}$, satisfying
$$F_{f_{1}*f_{2}}=F_{f_{2}}*F_{f_{1}},\quad f_{F_{1}*F_{2}}=f_{F_{2}}*f_{F_{1}}.$$ 
The harmonic analysis on commutative $I_{\pi_{j}}(G)$ can be further developed by introducing the spherical Fourier transform and the corresponding inversion formula cf.~\cite[\S3]{MR1450744}.

As we have observed before, the spherical functions $\Phi^{i}_{j}$ take values in the various algebras $\End(V_{i}(\pi_{j}))$. If $(G,H,\pi_{j})$ is a multiplicity free triple, then $V_{i}(\pi_{j})$ is $H$-isomorphic to $W_{j}$ if and only if $\dim\Hom_{H}(W_{j},V_{i})=1$. In this case it makes sense to denote
$$\widehat{G}(\pi_{j})=\{\rho_{i}\mid \dim\Hom_{H}(W_{j},V_{i})=1\}.$$
For $\rho_{i}\in\widehat{G}(\pi_{j})$ fix an $H$-equivariant embedding $J^{i}_{j}:W_{j}\to V_{i}$ and let $P^{i}_{j}:V_{i}\to W_{j}$ be its adjoint. Then
$$\Psi^{i}_{j}(x)=P^{i}_{j}\rho_{i}(x)J^{i}_{j}$$
is also called a spherical function of type $\pi_{j}$ associated to $\rho_{i}$ and it corresponds to the spherical function $\Phi^{i}_{j}$ of Definition \ref{Def:SF} upon the identification $W_{j}\cong V_{i}(\pi_{j})$ via $J^{i}_{j}$ and $P^{i}_{j}$.
The group $H\times H$ acts on the space of $\End(W_{j})$-valued functions by
$$(h_{1},h_{2})\Psi(x)=\pi_{j}(h_{1})\Psi(h_{1}^{-1}xh_{2})\pi_{j}(h_{2})^{-1}$$
so that the space of invariants,
$$E(\pi_{j})=\left(L(G)\otimes\End(W_{j})\right)^{H\times H},$$
consists of functions $\Psi:G\to\End(W_{j})$ that transform according to
$$\Psi(h_{1}xh_{2})=\pi_{j}(h_{1})\Psi(x)\pi_{j}(h_{2})$$
for all $h_{1},h_{2}\in H$ and $x\in G$. We have $\Psi^{i}_{j}\in E(\pi_{j})$ and following the Peter-Weyl decomposition $L(G)=\bigoplus_{j=1}^{m}\End(V_{j})$ we see that
$\{\Psi^{i}_{j}\mid\rho_{i}\in\widehat{G}(\pi_{j})\}$ is actually a basis of $E(\pi_{j})$.
Recall that $\pi_{1}$ is the trivial $H$-representation. The space $E(\pi_{1})$ is the algebra (for pointwise multiplication) of $H$-biinvariant functions on $G$ and $E(\pi_{i})$ is a module over $E(\pi_{1})$, for pointwise multiplication.
\begin{remark}\label{remark:notfree}
    We observe that if $E(\pi_{i})$ is freely and finitely generated as an $E(\pi_{1})$-module, then $\dim(E(\pi_{1}))$ divides $\dim(E(\pi_{i}))$.
\end{remark}


\section{The pair $(\GL_{2}(q)\times\GL_{2}(q),\mathrm{diag}(\GL_{2}(q)))$}\label{Section:product}

In this section we calculate the multiplicity table for the pair $(G\times G,\diag(G))$ with 
$G=\GL_{2}(q)$. First we recall the characters for $\GL_{2}(q)$ and their values on the conjugacy classes, where we follow \cite{Steinberg1951-Irreps}.
Then we proceed to decompose the tensor products and we present the results the multiplicity table. We defer most of the calculations to the appendix.

To discuss the conjugacy classes of $\GL_{2}(q)$ and the equivalence classes of irreducible representations of $\GL_{2}(q)$ we introduce the following notation. It will be convenient to work with the notation $r=q-1$ and $s=q+1$.

Consider the canonical projection $g:\bbZ_{rs}\to\bbZ_{s}$ and let $\calX'\subset\bbZ_{rs}$ be the complement of $\ker g$. Let $\bbZ_2$ act on both $\bbZ_{rs}$ and $\bbZ_{s}$ by letting the nontrivial element send $x\mapsto qx$. Since $g$ intertwines these actions and $\ker g$ is stable for this action, so is $\calX'$. The quotient $\calX'/\bbZ_{2}$ is denoted by $\calX$.
Note that $\calX$ has $(rs-r)/2=q(q-1)/2$ elements.

Let $\calW''=\bbZ_{r}\times\bbZ_{r}$ and consider the ring homomorphism $h:\calW''\to\bbZ_{r},(a,b)\mapsto a-b$. Let $\bbZ_{2}$ act on $\calW''$ by letting the non-trivial element send $(a,b)\mapsto(b,a)$. The complement $\calW'$ of $\ker h$ is stable for this action. The quotient $\calW'/\bbZ_{2}$ is denoted by $\calW$. Note that $\calW$ has $(r^2-r)/2=(q-1)(q-2)/2$ elements.

\subsection{Characters for $\GL_{2}(q)$.}
Let $\rho$ be a primitive element in $\bbF_{q}$ and $\sigma$ a primitive element in $\bbF_{q^{2}}$ with $\sigma^{s}=\rho$. The group $\GL_{2}(q)$ contains $q(q-1)^{2}(q+1)$ elements and 
four types of conjugacy classes, $c_{1}(k), c_{2}(k), c_{3}([k,\ell])$ or $c_{4}([n])$. Representatives, parameters, sizes and number of classes of these types are tabulated in Table \ref{GLcon}, where the representative of $c_{4}([n])$ is actually in $\GL(2,q^{2})$. 

\begin{table}[h]
\begin{center}
\begin{tabular}{|c c c c c|}
  \hline
  Conj.class & repr. & parameters & Size
& No. of classes \\ [0.5ex]
  \hline \hline
$c_1(k)$ & $\begin{pmatrix}
     \rho^k & 0 \\ 0 & \rho^k
\end{pmatrix}$ & $k\in\bbZ_{r}$ & $1$ & $r$\\
& & & & \\
$c_2(k)$ & $\begin{pmatrix}
     \rho^k & 0 \\ 1 & \rho^k
\end{pmatrix}$ & $k\in\bbZ_{r}$ & $rs$ & $r$\\
& & & & \\
$c_3([k,\ell])$ & $\begin{pmatrix}
     \rho^k & 0 \\ 0 & \rho^{\ell}
\end{pmatrix}$ & $[k,\ell]\in\calW$ & $qs$ & $\frac{r(r-1)}{2}$\\
& & & & \\
$c_4([n])$ & $\begin{pmatrix}
     \sigma^n & 0 \\ 0 & \sigma^{qn}
\end{pmatrix}$ & 

$[n]\in\calX$

& $qr$ & $\frac{qr}{2}$\\
\hline
  \end{tabular}
  \caption{Conjugacy classes of $\GL_{2}(q)$}\label{GLcon}
  \end{center}
\end{table}

The equivalence classes of irreducible representations of $\GL_{2}(q)$ are denoted by $U_{a}, V_{a}, W_{[b,c]}$ and $X_{[n]}$ whose parameters, dimensions and duals are tabulated in Table \ref{GLIrreps}.

\begin{table}[h]
    \centering
    \begin{tabular}{|c|c c c|}
    \hline
    Name & Dimension & Parameters & dual \\
    \hline
    $U_{a}$ & $1$ & $a\in\bbZ_{r}$ & $U_{-a}$\\
    $V_{a}$ & $q$ & $a\in\bbZ_{r}$ & $V_{-a}$\\
    $W_{[a,b]}$ & $s$ & $[a,b]\in\calW$ & $W_{[-a,-b]}$\\
    $X_{[n]}$ & $r$ & $[n]\in\calX$ & $X_{[-n]}$\\
    \hline
    \end{tabular}
    \caption{Equivalence classes of irreducible representations of $\GL_{2}(q)$.}
    \label{GLIrreps}
\end{table}

\begin{definition}
    For $a\in\bbZ_{r}$ the character $\alpha_{a}:\bbF_{q}^{\times}\to\bbC^{\times}$ is defined by $\alpha_{a}(\rho)=\exp(2a\pi i/r)$.
    For $n\in\bbZ_{rs}$ the character $\phi_{n}:\bbF^{\times}_{q^{2}}\to\bbC^{\times}$ is defined by
$\phi_{n}(\sigma)=\exp(2n \pi i/(rs))$.
\end{definition}
Note that $\alpha_{a}(\rho^{2j})=\alpha_{2a}(\rho^{j})$, $\alpha_{a}(\rho^{j})=\phi_{sa}(\sigma^{j})$ and $\phi_{n}(\sigma^{sj})=\alpha_{\overline{n}}(\rho^{j})$, where $\overline{n}$ is the image of $n\in\bbZ_{rs}$ under the natural projection $\bbZ_{rs}\to\bbZ_{r}$.
The values of the characters of $\GL_{2}(q)$ on the conjugacy classes are tabulated in Table \ref{GLchar}.

\begin{table}[h]
    \centering
    \begin{tabular}{|c|c c c c|}
    \hline
        $\GL_{2}(q)$ & $c_1(k)^{1}$ & $c_2(k)^{rs}$ & $c_3(k,\ell)^{qs}$ & $c_4([m])^{qr}$ \\
        \hline &&&&\\
        $U_a$ & $\alpha_a(\rho^{2k})$ & $\alpha_a(\rho^{2k})$ & $\alpha_a(\rho^{k+\ell})$ & $\alpha_a(\sigma^{sm})$ \\
        $V_a$ & $q\alpha_a(\rho^{2k})$ & $0$ & $\alpha_a(\rho^{k+\ell})$ & $-\alpha_a(\sigma^{sm})$ \\
        $W_{[a,b]}$ & $s\alpha_{a+b}(\rho^{k})$ & $\alpha_{a+b}(\rho^{k})$ & $\alpha_{a}(\rho^{k})\alpha_{b}(\rho^{\ell}) + \alpha_{a}(\rho^{\ell})\alpha_{b}(\rho^{k})$ & $0$ \\
        $X_{[n]}$ & $r\phi_{n}(\rho^k)$ & $-\phi_{n}(\rho^k)$ & $0$ & $-(\phi_{n}(\sigma^m) + \phi_{n}(\sigma^{qm}))$ \\
        &&&&\\
        \hline
    \end{tabular}
    \caption{Character table for $\GL_{2}(q)$. The exponents of the conjugacy classes indicate the number of elements in this class.}
    \label{GLchar}
\end{table}

\subsection{Tensor product decomposition}

The irreducible representations of $G\times G$ are of the form $\pi_{1}\otimes\pi_{2}$ with $\pi_{1},\pi_{2}$ irreducible representations of $G$. Note that the restriction of $\pi_{1}\otimes\pi_{2}$ and $\pi_{2}\otimes\pi_{1}$ to the diagonal $\diag(G)$ are the same. 

To calculate the multiplicities we identify $\diag(G)$ with $G$ itself and we consider an irreducible representation of $G\times G$ as a $G$-representation in the tensor product of two $G$-representations. The character of such a representation $\pi_{1}\otimes\pi_{2}$ is just the product of the characters, $\chi_{\pi_{1}}\chi_{\pi_{2}}$. To calculate the multiplicity of $\pi_{3}$ inside $\pi_{1}\otimes\pi_{2}$ we invoke
\eqref{FormulaMultiplicity}. However, we first discuss some observations.

\begin{example}\label{ExampleTensorU}
    The tensor product of $U_{a}$ with any irreducible representation remains irreducible. To determine the irreducible representation that we obtain we compare the characters directly, see Table \ref{TensorWithU}.  

    \begin{table}[ht]
    \centering
    \begin{tabular}{c|c c c c|}
        
        & $c_1(k)$ & $c_2(k)$ & $c_3([k,\ell])$ & $c_4([m])$ \\
        \hline &&&&\\
        $U_a \otimes U_b$ & $\alpha_{a+b}(\rho^{2k})$ & $\alpha_{a+b}(\rho^{2k})$ & $\alpha_{a+b}(\rho^{k+\ell})$ & $\alpha_{a+b}(\sigma^{sm})$ \\
        $U_a \otimes V_b$ & $q\alpha_{a+b}(\rho^{2k})$ & $0$ & $\alpha_{a+b}(\rho^{k+\ell})$ & $-\alpha_{a+b}(\sigma^{sm})$ \\
        &&&& \\
        $U_a \otimes W_{[b,c]}$ & $s\alpha_{2a+b+c}(\rho^{k})$ & $\alpha_{2a+b+c}(\rho^{k})$ & $\makecell{\alpha_{a+b}(\rho^{k})\alpha_{a+c}(\rho^{\ell}) ~+ \\ \alpha_{a+b}(\rho^{\ell})\alpha_{a+c}(\rho^{k})}$ & $0$ \\
        &&&& \\
        $U_a \otimes X_{[n]}$ & $r\phi_{n+sa}(\rho^k)$ & $-\phi_{n+sa}(\rho^k)$ & $0$ & $\makecell{-(\phi_{n+sa}(\sigma^m)~ + \\ \phi_{n+sa}(\sigma^{qm}))}$ \\
        &&&&\\
        \hline
    \end{tabular}
    \caption{The characters of $U_{a}\otimes\pi$ for $\pi$ irreducible.}\label{TensorWithU}
\end{table}
Based on Table \ref{TensorWithU} we observe the following equalities:

\[
U_a \otimes U_b = U_{a + b}, \quad
U_a \otimes V_b = V_{a + b}, \quad
U_a \otimes W_{[b,c]} = W_{[a + b, \, a + c]}, \quad
U_a \otimes X_{[n]} = X_{[n + sa]}.
\]
\end{example}

\begin{lemma}\label{sym}
For irreducible representations $\pi_{1},\pi_{2},\pi_{3}$ we have
$$\dim\Hom_{G}(V_{\pi_{1}}\otimes V_{\pi_{2}},V_{\pi_{3}})=\dim\Hom_{G}(V_{\pi_{1}}\otimes V_{\pi_{3}}^{*},V_{\pi_{2}}^{*}).$$
\end{lemma}
\begin{proof}
    This follows from the canonical isomorphism $(V_{\pi_{1}}\otimes V_{\pi_{2}})^{*}\otimes V_{\pi_{3}}=(V_{\pi_{1}}\otimes V_{\pi_{3}}^{*})^{*}\otimes V_{\pi_{2}}^{*}$ of $G$-modules and taking $G$-invariants.
\end{proof}
\begin{example}\label{ExampleVoV}
    We have $\dim\Hom_{G}(V_{a}\otimes V_{b},U_{b'})=\dim\Hom_{G}(V_{a}\otimes U_{-b'},V_{-b})$ which is less than or equal to one with equality if and only if $a-b'=-b$. It follows that the only constituent of type $U_{b'}$ in the $G$-module $V_{a}\otimes V_{b}$ is $U_{a+b}$ and the multiplicity is one.
\end{example}

In general we need \eqref{FormulaInnerProduct} and \eqref{FormulaMultiplicity} to calculate the multiplicities but instead of summing over the elements in $G$ we can also sum over the conjugacy classes,
\begin{equation}\label{MultFormulaB}
    \dim\Hom_{G}(V_{\pi_{1}}\otimes V_{\pi_{2}},V_{\pi_{3}})=\frac{1}{|G|}\sum_{c\in G/\Ad(G)}|c|\chi_{\pi_{1}}(c)\chi_{\pi_{2}}(c)\overline{\chi_{\pi_{3}}(c)}
\end{equation}
where $\chi(c)$ is just the evaluation of $\chi$ on an element of the conjugacy class $c$.

\begin{lemma}\label{GLformula}
The multiplicity of $\pi_{3}$ in the tensor product $\pi_{1}\otimes\pi_{2}$ is given by
\[
\begin{aligned}
\dim \Hom_G\big(V_{\pi_1} \otimes V_{\pi_2}, \, V_{\pi_3}\big)
&= \frac{1}{q s r^2} \bigg( 
\sum_{k\in\bbZ_{r}} \chi_{\pi_1}(c_1(k)) \, \chi_{\pi_2}(c_1(k)) \, \overline{\chi_{\pi_3}(c_1(k))} \\
&\quad + \, rs \sum_{k\in\bbZ_{r}} \chi_{\pi_1}(c_2(k)) \, \chi_{\pi_2}(c_2(k)) \, \overline{\chi_{\pi_3}(c_2(k))} \\
&\quad + \, qs \sum_{[k,\ell]\in\calW} \chi_{\pi_1}(c_3([k,\ell])) \, \chi_{\pi_2}(c_3([k,\ell])) \, \overline{\chi_{\pi_3}(c_3([k,\ell]))} \\
&\quad + \, {qr} \sum_{[n]\in\calX} \chi_{\pi_1}(c_4([n])) \, \chi_{\pi_2}(c_4([n])) \, \overline{\chi_{\pi_3}(c_4([n]))}
\bigg).
\end{aligned}
\]
\end{lemma}
\begin{proof}
    This follows directly from \eqref{MultFormulaB} and the data in Table \ref{GLcon}.
\end{proof}

\begin{remark}
    From now we denote for three irreducible representations $\pi_{1},\pi_{2},\pi_{3}$,
    $$[\pi_{1}\otimes\pi_{2},\pi_{3}]=\dim \Hom_G\big(V_{\pi_1} \otimes V_{\pi_2}, \, V_{\pi_3}\big).$$
    To calculate the decomposition of the tensor product into irreducible submodules, we have to determine the multiplicities $[\pi_{1}\otimes\pi_{2}:\pi_{3}]$ where $\pi_{3}$ ranges over the irreducible $G$ representations $U_{a},V_{b},W_{[c,d]}$ and $X_{[n]}$ for the various parameters.
\end{remark}

\begin{remark}\label{RemarkCalc}
    In the calculations that follow we will use some identities of sums of characters. The indicator function $\mathbbm{1}$ is one if the condition is satisfied and zero else. 
    Orthogonality of characters implies
    $$\sum_{j \in \mathbb{Z}_{r}} \alpha_a(\rho^j)= r \, \mathbbm{1}_{a = 0},\quad \sum_{j \in \mathbb{Z}_{rs}} \phi_{n}(\sigma^j)= rs \, \mathbbm{1}_{n = 0}.$$
    Almost direct applications of these formulas yield the identities
    $$\sum_{[j] \in \mathcal{X}} (\phi_n(\sigma^j)
    + \phi_n(\sigma^{qj}))
    = rs \, \mathbbm{1}_{n = 0} - r \, \mathbbm{1}_{\overline{n} = 0},
    \quad \text{where } \overline{n} \in \mathbb{Z}_{r}$$
    and
    $$\sum_{[j] \in \mathcal{X}} \alpha_a(\sigma^{s j})
    = \frac{1}{2} \left( rs \, \mathbbm{1}_{a = 0} - r \, \mathbbm{1}_{2a = 0} \right).$$
    Finally,
    $$\sum_{[j,k]\in\calW}\left(\alpha_{a}(\rho^{j})\alpha_{b}(\rho^{k})+\alpha_{a}(\rho^{k})\alpha_{b}(\rho^{j})\right)
    =r^2\mathbbm{1}_{a=0,b=0}-r\mathbbm{1}_{a+b=0}$$
    and in particular
    $$\sum_{[j,k] \in \mathcal{W}} \alpha_a(\rho^{j+k})
    = \frac{1}{2} \left( r^2 \, \mathbbm{1}_{a = 0} - r \, \mathbbm{1}_{2a = 0} \right).$$
\end{remark}
With these identities  we can calculate the other multiplicities. In some calculations we do not use these identities, but rather the ideas of their proof. For example, if we sum over elements in $\calX$, then instead we sum over $\bbZ_{rs}$, subtract what we have counted too much, and divide by two to account for the equivalence relation.

\begin{lemma}\label{Lemma:VV}
    The decomposition of $V_{a}\otimes V_{b}$ is described by
    \begin{itemize}
        \item $[V_{a}\otimes V_{b}:U_{a'}]=\mathbbm{1}_{a'=a+b}$,
        \item $[V_{a}\otimes V_{b}:V_{b'}]=\mathbbm{1}_{2b'=2(a+b)}$,
        \item $[V_{a}\otimes V_{b}:W_{[c',d']}]=\mathbbm{1}_{c'+d'=2(a+b)}$
        \item $[V_{a}\otimes V_{b}:X_{[n']}]=\mathbbm{1}_{\overline{n'}=2(a+b)}$.
    \end{itemize}
\end{lemma}

\begin{proof}
    For $[V_{a}\otimes V_{b}:U_{a'}]$ see Example \ref{ExampleVoV}. For $[V_{a}\otimes V_{b}:V_{b'}]$ we use Lemma \ref{GLformula} and the data Table \ref{GLchar},
    \[
    \begin{aligned}
    qsr^2[V_a \otimes V_b : V_{b'}] = \; & q^3 \sum_{j\in\bbZ_{r}} \alpha_{a+b-b'}(\rho^{2j}) + qs \sum_{[j,k]\in\calW} \alpha_{a+b-b'}(\rho^{j+k})- qr \sum_{[j]\in\calX} \alpha_{a+b-b'}(\sigma^{sj}) \\
    =& q^3 r \mathbbm{1}_{2a+2b=2b'}+\frac{1}{2}qs(r^2\mathbbm{1}_{a+b=b'}-r\mathbbm{1}_{2a+2b=2b'})-\frac{1}{2}qr(rs\mathbbm{1}_{a+b=b'}-r\mathbbm{1}_{2a+2b=2b'})\\
    =&qsr^2\mathbbm{1}_{2a+2b=2b'}
    \end{aligned}
    \]    
    We continue with $[V_a \otimes V_b : W_{[c',d']}]$,
    $$
    \begin{aligned}
    q s r^2 \, [V_a \otimes V_b : W_{[c',d']}] = & \, q^2 s \sum_{j\in\bbZ_{r}} \alpha_{2a+2b-c'-d'}(\rho^{j}) \\
    & + q s \sum_{[j,k]\in\calW} \alpha_a(\rho^{j+k}) \alpha_b(\rho^{j+k}) \left( \alpha_{-c'}(\rho^j) \alpha_{-d'}(\rho^k) + \alpha_{-c'}(\rho^k) \alpha_{-d'}(\rho^j) \right) \\
    = & \, q^2 s r \mathbbm{1}_{2a+2b=c'+d'} \\
    & + q s  \sum_{[j,k]\in\calW}\left(\alpha_{a + b - c'}(\rho^j) \alpha_{a + b - d'}(\rho^k) +\alpha_{a + b - c'}(\rho^k) \alpha_{a + b - d'}(\rho^j)  \right) \\
    =&  (q^2 s r -qsr)\mathbbm{1}_{2a+2b=c'+d'} = q s r^2 \mathbbm{1}_{2a+2b=c'+d'}.
    \end{aligned}
    $$
    Finally,
    $$
    \begin{aligned}
    qsr^2 [V_a \otimes V_b : X_{[n']}] = & \, rq^2 \sum_{j\in\bbZ_{r}} \alpha_{2a+2b-\overline{n'}}(\rho^{j})\\
    & - qr \sum_{[j]\in\calX} \alpha_a(\sigma^{sj}) \alpha_b(\sigma^{sj}) \left( \phi_{-n'}(\sigma^j) + \phi_{-n'}(\sigma^{qj}) \right) \\
    = & \, q^2r^2 \mathbbm{1}_{2a + 2b =\overline{ n'}} -qr \sum_{[j]\in\calX}(\phi_{sa + sb - n'}(\sigma^j)+\phi_{sa + sb - n'}(\sigma^{qj})) \\
    = & \, (q^2 r^2\mathbbm{1}_{2a + 2b = \overline{n'}} -qr(rs\mathbbm{1}_{sa+sb=n'}-r\mathbbm{1}_{sa+sb=\overline{n'}})\\
    = & \, q s r^2 \, (\mathbbm{1}_{2a + 2b = \overline{n'}}-\mathbbm{1}_{sa+sb=n'})=q s r^2 \, \mathbbm{1}_{2a + 2b = \overline{n'}},
    \end{aligned}
    $$
    because $n'=s(a+b)$ contradicts $[n']\in\calX$.
\end{proof}

\begin{proposition}
    The multiplicity table for the pair $(G\times G,\diag(G))$ is given in Table \ref{table:tensorproductdecomp}.
\end{proposition}

\begin{proof}
    The first five rows of Table \ref{table:tensorproductdecomp} follow from Example \ref{ExampleTensorU} and Lemma \ref{Lemma:VV}. All other rows follow from the results in Appendix \ref{Appendix:GxG}. 
\end{proof}

\begin{remark}\label{TableExplanation}
    The entries in Table \ref{table:tensorproductdecomp} are equations that represent the conditions of the indicator functions. Simultaneous conditions are separated by a comma, the plus indicates a sum of indicator functions.
    The equations are expressed in the parameters of the classes of irreducible representations. If all parameters are in $\bbZ_{r}$ then the equation is in $\bbZ_{r}$, for example $a+b=b'$ or $2(a+b)=\overline{n'}$ (recall that $\overline{n'}\in\bbZ_{r}$). Some equations are in $\bbZ_{rs}$, for example $sa+n=n'$. 
\end{remark}

\begin{table}
\centering
\scriptsize
\begin{tabular}{|c|c|c|c|c|}
\hline
 & $U_{a'}$ & $V_{b'}$ & $W_{[c',d']}$ & $X_{[n']}$  \\
\hline
\hline
$U_{a}\otimes U_{b}$ & $a+b=a'$ & $0$ & $0$ & $0$ \\
\hline
$U_{a}\otimes V_{b}$ & $0$ & $a+b=b'$ & $0$ & $0$ \\
\hline
$U_{a}\otimes W_{[b,c]}$ & $0$ & $0$ &

\(\begin{array}{c}
a+b=c',a+c=d' \\
+ \\
a+c = c',a+b = d'
\end{array}\) 

& $0$ \\
\hline
$U_{a}\otimes X_{[n]}$ & $0$ & $0$ & $0$ & $\begin{array}{c}
s a + n = n' \\
\end{array}$ \\
\hline
$V_{a}\otimes V_{b}$ & $a+b=a'$ & $2(a+b) = 2b'$ & $2(a+b) = c'+ d'$ & $2(a+b) = \overline{n'}$ \\
\hline
$V_a \otimes W_{[b,c]} $ & $0$ & $ 2a +b +c = 2b'$ &  \(\begin{array}{c}
2a+b+c=c'+d' \\
+ \\
a+b = c',a+c = d' \\
+ \\
a+c = c',a+b = d'
\end{array}\)  &  $ 2a + b + c = \overline{n'}$ \\

\hline
$V_a \otimes X_{[n]}$ & $0$ & $ 2a + \overline{n} = 2b'$ & $2a + \overline{n} = c'+ d'$ &  \(\begin{array}{c}
2a + \overline{n} = \overline{n'} \\
- \\ 
sa + n = n '  \\
- \\
sa + n = q n '
\end{array}\) \\

\hline
\(
W_{[a,b]} \otimes W_{[c,d]}\) 
& $\begin{array}{c}
a+c=a',b+d=a'\\
+\\
a+d=a',b+c=a'
\end{array}$ & \(\begin{array}{c}
2b'=a+b+c+d\\+\\b'=a+c,b'=b+d \\
+ \\
b'=a+d,b'=b+c
\end{array}\) &
\(\begin{array}{c}
a+b+c+d = c'+d' \\
+ \\
a+c = c', b+d=d' \\
+ \\
a+c=d' , b+d = c' \\
+ \\
a+d = c', b+c=d' \\
+ \\
a+d=d' , b+c = c'
\end{array}\)
& $a+b+c+d = \overline{n'}$ \\
\hline
$W_{[a,b]} \otimes X_{[n]}$ & $0$ & $a+b+\overline{n} = 2b'$ & $a+b+\overline{n} = c'+ d'$ & $a+b+\overline{n} = \overline{n'}$ \\
\hline
$X_{[n]} \otimes X_{[m]}$
& $\begin{array}{c}
n + m = sa'
\end{array}$  & 
$\begin{array}{c}
\overline{n} + \overline{m} = 2b' \\
- \\
n + m = sb' \\
- \\
qn + m = sb'
\end{array}$
 & $\overline{n} + \overline{m} = c'+ d'$& \(\begin{array}{c}
\overline{n} + \overline{m} = \overline{n'} \\
- \\
n + m = n '  \\
- \\
qn + m = n ' \\
- \\
n + qm = n '  \\
- \\
n + m = qn ' 
\end{array}\) \\
\hline
\end{tabular}
\caption{The multiplicity table of $\GL_{2}(q)\times\GL_{2}(q),\mathrm{diag}(\GL_{2}(q))$, see Remark \ref{TableExplanation} for the explanation of the entries}
\label{table:tensorproductdecomp}
\end{table}

\subsection{Application}

As an application of the multiplicities in Table \ref{table:tensorproductdecomp} we identify the Gelfand triples for the pair $(\GL_{2}(q)\times\GL_{2}(q),\diag(\GL_{2}(q))$.

\begin{remark}\label{remark:NoGelfandTripleVW}
    We have $[V_{a}\otimes W_{[b,c]}:W_{[a+b,a+c]}]=2$,
    which means that no tensor product $V_{a}\otimes W_{[b,c]}$ decomposes multiplicity free. Conversely, no module $W_{[b,c]}$ induces multiplicity free because $[V_{0}\otimes W_{[b,c]}:W_{[b,c]}]=2$. Likewise, no module $V_{a}$ induces multiplicity free.
\end{remark}

\begin{proposition}
    Let $\pi$ be an irreducible representation of $\GL_{2}(q)$ that we view as an irreducible representation of $\diag(\GL_{2}(q))\subset\GL_{2}(q)\times\GL_{2}(q)$. Then
    $(\GL_{2}(q)\times\GL_{2}(q),\diag(\GL_{2}(q)),\pi)$ is a Gelfand triple if and only if $\dim(\pi)=1$ or $\dim(\pi)=q-1$.
\end{proposition}

\begin{proof}
    If $\dim(\pi)=1$ of $\dim(\pi)=q-1$ then $\pi$ occurs in the decomposition of each irreducible representation of $\GL_{2}(q)\times\GL_{2}(q)$ with multiplicity at most one, which follows from Table \ref{table:tensorproductdecomp}. Indeed, in the columns of $U_{a'}$ and $X_{n'}$ there are sums of indicator functions in some entries, but the equations of the indicator functions cannot hold simultaneously. 
    
    Conversely, for each representation of dimension $q$ or $q+1$ there is an irreducible representation of $\GL_{2}(q)\times\GL_{2}(q)$ that contains this irreducible representation with multiplicity two, as we have observed in Remark \ref{remark:NoGelfandTripleVW}.
\end{proof}

The $\GL_{2}(q)\times\GL_{2}(q)$ module $\Ind U_{0}$ decomposes into $q^2-1$ irreducible submodules $\pi^{*}\otimes\pi$ with $\pi\in\widehat{G}$, its dimension is $\dim\Ind U_0=|\GL_{2}(q)|=qsr^2$. In general we have
$$\Ind U_{a'}=\bigoplus_{a\in\bbZ_{r}}U_{a}\otimes U_{a'-a}\oplus\bigoplus_{a\in\bbZ_{r}}V_{a}\otimes V_{a'-a}\oplus\bigoplus_{(a,b)\in\calW}W_{[a,b]}\otimes W_{[a'-a,a'-b]}\oplus\bigoplus_{[n]\in\calX}X_{[n]}\otimes X_{[sa'-n]},$$
a module of dimension $qsr^2$ that decomposes multiplicity free into $q^2-1$ irreducible submodules.

\begin{remark}
It follows that the $\bbC$-algebra of zonal spherical functions $E(U_{0})$ is a complex vector space of dimension $q^2-1$. The $E(U_{0})$-module $E(U_{a})$ is of the same dimension, but its module structure not clear at the moment.
\end{remark}

Concerning the decomposition of $\Ind X_{[n']}$ we are interested in the number of irreducible constituents and also in the constituents themselves. The number of constituents is a counting problem and the results are displayed in Table \ref{TableIndX}. Most of the items in the table are easy to reproduce using Table \ref{table:tensorproductdecomp}, but for a couple of items the counting is a bit involved. We discuss a couple of the more involved entries below.

\begin{example}\label{exampleWW}
    We will first count the number of representations of type $W_{[a,b]} \otimes W_{[c,d]}$ occuring in the decomposition of $\Ind X_{[n']}$, or the number of pairs $([a,b],[c,d]) \in \calW \times \calW$ for which $X_{[n']}$ is a submodule of $W_{[a,b]} \otimes W_{[c,d]}$. \\

    We first consider the case that $p$ is odd and $n'$ even. If $a+b$ is odd, then $c+d$ has to be odd as well; in particular, we cannot have $c=d$, so we do not have to worry about reducible representations. We have $\frac{1}{2}r^2$ choices for $a$ and $b$ such that $a+b$ is odd, $r$ choices for $a$ and $\frac{1}{2}r$ choices for $b$ since it must be of a different parity than $a$, and then $r$ choices for $c$: this fixes the choice for $d$, giving us $\frac{1}{2}r^3$ such modules.
    
    We now count the number of pairs for which $a+b$ is even: we have $r(\frac{1}{2}r-1) $ choices for such $a$ and $b$; $r$ choices for $a$ and then $\frac{1}{2}r-1$ choices for $b$, as it must be of the same parity as $a$, but not equal to $a$. Since $a+b$ is even, so is $c+d$. We now have $r-2$ choices for $c$, as exactly $2$ choices will force $d$ to be equal to $c$, giving a reducible representation, giving us another $r(\frac{1}{2}r-1)(r-2)$ modules.
    
    Finally, we divide the total number of such modules by four, as we have counted every module four times by making no distinction between $((a,b),(c,d)) $, $((b,a),(c,d))$, $((a,b),(d,c))$ and $((b,a),(d,c))$.
    We conclude that the number of pairs $([a,b],[c,d]) \in \calW \times \calW$ for which $X_{[n']}$ is a submodule of $W_{[a,b]} \otimes W_{[c,d]}$ with $n'$ even is equal to
    $$\frac{1}{4}(\frac{1}{2}r^3 + r(\frac{1}{2}r-1)(r-2)) = \frac{1}{4}(q-1)(q^2-4q+5).$$

    The count for the case where $p$ and $n'$ are odd is similar: the difference with the calculation above lies in the following: when $a+b$ is odd, $c+d$ is even, giving us $\frac{1}{2}r^2$ choices for $(a,b)$ and $r-2$ choices for $c$, as $d$ is fixed and $2$ choices of $c$ result in $c=d$. When $a+b$ is even, $c+d$ is odd, giving us $r(\frac{1}{2}r-1)$ choices for $(a,b)$ and $r$ choices for $c$. After dividing by four, we conclude that the number of pairs $([a,b],[c,d]) \in \calW \times \calW$ for which $X_{[n']}$ is a submodule of $W_{[a,b]} \otimes W_{[c,d]}$ with $n'$ odd is equal to
    $$\frac{1}{4}(\frac{1}{2}r^2(r-2) +r(\frac{1}{2}r-1)r ) = \frac{1}{4}(q-1)(q^2-4q+3).$$

    Finally, we consider the case where $p$ is even. Now we no longer have a distinction for $a+b$ is odd or even. We have $r$ choices for $a$, $r-1$ choices for $b$ such that $a \neq b$ and $r-1$ choices for $c$ such that $c \neq d$. After dividing by four, we conclude that for $p$ even, the number of pairs $([a,b],[c,d]) \in \calW \times \calW$ for which $X_{[n']}$ is a submodule of $W_{[a,b]} \otimes W_{[c,d]}$ is equal to $$\frac{1}{4}r(r-1)^2 = \frac{1}{4}(q-1)(q-2)^2.$$
\end{example}

\begin{example}\label{exampleWX}
    We will now count the number of pairs $([a,b],[n]) \in \calW \times \calX$ for which $X_{[n']}$ is a submodule of $W_{[a,b]} \otimes X_{[n]}$ and $X_{[n]} \otimes W_{[a,b]}$. We first consider the case $p$ is odd and $n'$ is even. First take $a+b$ to be odd: this gives us $r$ choices for $a$ and $\frac{1}{2}r$ choices for $b$, since $b$ must be of a different parity than $a$. Since $n$ now has to be odd, we cannot have $n = 0$ modulo $s$, thus we have $s$ choices for $n$ such that $a+b+n = n '$ modulo $r$, giving us a total of $\frac{1}{2}r^2s$ pairs. Now take $a+b$ to be even: we then have $r$ choices for $a$ and $\frac{1}{2}r-1$ choices for $b$, as $b$ has to be of the same parity as $a$, but we need $b \neq a$. Since $n$ now has to be even, we have to exclude the choices for $n$ that satisfy $n = 0$ modulo $s$, giving us $s-2$ choices of $n$ such that
    $a+b+n = n'$ modulo $r$, giving us another $$r(\frac{1}{2}r-1)(s-2)$$ pairs. 
    Lastly, we divide the number by four, as we have counted every module four times by making no distinction between $((a,b),n) $, $((b,a),n)$, $((a,b),qn)$ and $((b,a),qn)$.
    We conclude that the number of pairs $([a,b],n) \in \calW \times \calX$ for which $X_{[n']}$ is a submodule of $W_{[a,b]} \otimes X_{[n]}$ for $n'$ even is equal to $$\frac{1}{4}(\frac{1}{2}r^2s + r(\frac{1}{2}r-1)(s-2)) = \frac{1}{4}(q-1)^3 .$$
    We get the same number for the modules $X_{[n]} \otimes W_{[a.b]}$.

    Now consider the case where $p$ and $n'$ are odd. If $a+b$ is even, which is the case for $r(\frac{1}{2}r-1)$ choices of $(a,b)$, then $n$ is odd and in particular $n \neq 0$ modulo $s$, resulting in $s$ choices for $n$ such that $a+b+n = n'$ modulo $r$, giving us $r(\frac{1}{2}r-1)s$ pairs. If $a+b$ is odd, which is the case for $\frac{1}{2}r^2$ choices of $(a,b)$, we have to exclude the cases $n = 0$ modulo $s$, giving us $s-2$ choices for $n$ such that $a+b+n=n'$ modulo $r$, giving us another $\frac{1}{2}r^2(s-2)$ pairs. After diving the total number by four, we conclude that the number of pairs $([a,b],[n]) \in \calW \times \calX$ for which $X_{[n']}$ is a submodule of $W_{[a,b]} \otimes X_{[n]}$ for $n'$ odd is equal to $$\frac{1}{4}(r(\frac{1}{2}r-1)s + \frac{1}{2}r^2(s-2)) = \frac{1}{4}(q-1)(q^2-2q-1) .$$
    We get the same number for the modules $X_{[n]} \otimes W_{[a,b]}$.

    Finally we consider the case where $p$ is even. This gives us $r$ choices for $a$, $r-1$ choices for $b$ since we need $b \neq a$ and $s-1$ choices for $n$ such that $a+b+n=n'$ modulo $r$ since we need $n \neq 0$ modulo $s$. After diving by four, we conclude that the number of pairs $([a,b],n) \in \calW \times \calX$ for which $X_{[n']}$ is a submodule of $W_{[a,b]} \otimes X_{[n]}$ for $p$ even is equal to $$\frac{1}{4}r(r-1)(s-1) = \frac{1}{4}q(q-1)(q-2).$$
    We get the same number for the modules $X_{[n]} \otimes W_{[a,b]}$.
\end{example}

\begin{example}\label{exampleXX}  
We want to find the number of modules of the form $X_{[m]}\otimes X_{[n]}$ for which $[X_{[m]}\otimes X_{[n]}:X_{[n']}]=1$. We first consider the case that $p$ is odd and $n'$ is even.

If $m$ is odd ($sr/2$ possibilities), then $n$ has to be odd and therefore we do not have to worry about $n = 0$ modulo $s$. Then there are $s$ choices of $n\in\bbZ_{rs}$ for which $\overline{m}+\overline{n}=\overline{n'}$, but using Table \ref{table:tensorproductdecomp} we see that we have four exceptions. Hence there are $\tfrac12sr(s-4)$ such pairs.

If $m$ is even ($\tfrac12rs-r$ choices), then $n$ has to be even as well. The number of such $n\in\bbZ_{rs}$, not a multiple of $s$, for which $\overline{m}+\overline{n}=\overline{n'}$, depends on $m$:
$$\begin{array}{rcl}
s-2 & \mbox{if} & n'+m= n'-m=0\pmod{s},\\
s-4 & \mbox{if} & n'+m= 0 \pmod{s} ~~\text{or}~~n'-m=0\pmod{s},\\
s-6 & \mbox{else.} &
\end{array}$$
The first case only holds for $r$ choices of $m$ when $n' = \tfrac{1}{2}s$ modulo $s$, assuming $\tfrac{1}{2}s$ is even.
The second case holds for $2r$ choices of even $m$, assuming $n' \neq \frac{1}{2}s$ modulo $s$, and the third holds for all other even $m$; $\tfrac12rs-2r$ choices if $n = \frac{1}{2}s$ modulo $s$ and $\frac{1}{2}s$ is even and $\tfrac12rs-3r$ choices if $n \neq \frac{1}{2}s$ modulo $s$. This makes a total of
$$2r(s-4)+(\tfrac12rs-3r)(s-6)$$
of such pairs if $n' \neq \frac{1}{2}s$ modulo $s$. If $n' = \frac{1}{2}s$ modulo $s$ and $\frac{1}{2}s$ is even, then we have $$r(s-2)+(\tfrac12rs-2r)(s-6)$$ of such pairs, and we conclude that the number of pairs is the same for all $n'$ even.

Finally, we divide the number of pairs by four as we have counted every module four times by making no distinction between $(n,m), ~(qn,m),~(n,qm) $ and $(qn,qm)$.
We conclude that the number of pairs $([m],[n])\in\calX\times\calX$ for which $X_{[n']}$ is a submodule of $X_{[m]}\otimes X_{[n]}$ for $n'$ even is equal to $$\frac{1}{4}((\frac{1}{2}sr+2r)(s-4) + (\frac{1}{2}sr-3r)(s-6))  =\frac{1}{4}(q-1)(q^2-4q+5).$$

If $p$ and $n'$ are odd, then a similar count can be made. Now if $m$ is even, which gives us $\frac{1}{2}sr-r$ possibilities, $n$ has to be odd and we do not have to worry about $n = 0$ modulo $2$. We have $s-4$ choices for $n$, using the exceptions in Table \ref{table:tensorproductdecomp}. If $m$ is odd, which gives us $\frac{1}{2}rs$ choices, $n$ must be even. The number of such $n$, not a multiple of $s$, for which $\overline{m}+\overline{n}=\overline{n'}$, depends on $m$:
$$\begin{array}{rcl}
s-2 & \mbox{if} & n'+m= n'-m=0\pmod{s},\\
s-4 & \mbox{if} & n'+m= 0 \pmod{s} ~~\text{or}~~n'-m=0\pmod{s},\\
s-6 & \mbox{else.} &
\end{array}$$
The first case only holds for $r$ choices of $m$ when $n' = \tfrac{1}{2}s$ modulo $s$, assuming $\tfrac{1}{2}s$ is odd.
The second case holds for $2r$ choices of odd $m$ when $n' \neq \tfrac{1}{2}s$ modulo $s$, and the third holds for all other odd $m$; $\tfrac12rs-r$ choices if $n = \frac{1}{2}s$ modulo $s$ and $\frac{1}{2}s$ is odd and $\tfrac12rs-2r$ choices if $n \neq \frac{1}{2}s$ modulo $s$. This makes a total of
$$2r(s-4)+(\tfrac12rs-2r)(s-6)$$
of such pairs if $n' \neq \frac{1}{2}s$ modulo $s$. If $n' = \frac{1}{2}s$ modulo $s$ and $\frac{1}{2}s$ is odd, then we have $$r(s-2)+(\tfrac12rs-r)(s-6)$$ of such pairs, and we conclude that the number of pairs is the same for all $n'$ odd. After dividing by four, we conclude that the number of pairs $([m],[n])\in\calX\times\calX$ for which $X_{[n']}$ is a submodule of $X_{[m]}\otimes X_{[n]}$ for $n'$ odd is equal to $$\frac{1}{4}((\frac{1}{2}sr+r)(s-4)+(\frac{1}{2}sr-2r)(s-6))=\frac{1}{4}(q-1)(q^2-4q+3).$$

Finally, we consider the case where $p$ is even. The number of choices we have for $n$ depends on the choice of $m$:
$$\begin{array}{rcl}
s-3 & \mbox{if} & n'+m= 0 \pmod{s} ~~\text{or}~~n'-m=0\pmod{s},\\
s-5 & \mbox{else.} &
\end{array}$$
The first case holds for $2r$ choices of $m$, and the second case for all others. After dividing by four, we conclude that the number of pairs $([m],[n])\in\calX\times\calX$ for which $X_{[n']}$ is a submodule of $X_{[m]}\otimes X_{[n]}$ for $p$ even is equal to $$\frac{1}{4}(2r(s-3) + r(r-1)(s-5))  = \frac{1}{4}(q-1)(q-2)^2.$$

\end{example}

\begin{table}[h]
    \centering
    \begin{tabular}{c|cc|c}
        & \multicolumn{2}{c|}{$q$ odd} & $q$ even \\
        dimension & $n'$ even & $n'$ odd & \\
        \hline
        $r$ & $2(q-1)$ & $2(q-1)$ & $2(q-1)$\\
        $q^2$ & $2(q-1)$ & $0$ & $(q-1)$ \\
        $qs$ & $(q-1)(q-3)$ & $(q-1)^2$ & $(q-1)(q-2)$\\
        $qr$ & $(q-1)(q-3)$ & $(q-1)^2$ & $(q-1)(q-2)$\\
        $s^2$ & $\frac{1}{4} (q-1)(q^2-4q+5)$ & $\frac{1}{4} (q-1)(q^2-4q+3)$ & $\frac{1}{4} (q-1)(q-2)^2$ \\
        $rs$ & $\frac{1}{2}(q-1)^3$ & $\frac{1}{2}(q-1)(q^2-2q-1)$ & $\frac{1}{2}q(q-1)(q-2)$ \\
        $r^2$ & $\frac{1}{4} (q-1)(q^2-4q+5)$ & $\frac{1}{4} (q-1)(q^2-4q+3)$ & $\frac{1}{4} (q-1)(q-2)^2$
    \end{tabular}
    \caption{The number of irreducible submodules of $\Ind X_{[n']}$ of given dimension. For $q$ odd, the cases $n'$ even and $n'$ odd are distinguished.}
    \label{TableIndX}
\end{table}

We conclude with two observations concerning the entries of Table \ref{TableIndX}.
\begin{itemize}
    \item The number of tensor products of irreducible representations that contain $X_{[n']}$ is equal to $(q-1)(q^2-q+1)$. However, the module $\Ind X_{[n']}$ does depend on $q$ and on the parity of $n'$ in the case that $q$ is even.
    \item The $E(U_{0})$-module $E(X_{[n']})$ of spherical functions of type $X_{[n']}$, is a complex vector space of dimension $(q-1)(q^2-q+1)$, but it is not a free module because $\dim E(U_{0})=q^2-1$, see Remark \ref{remark:notfree}.
\end{itemize}


\section{The pair $(\SL_{q}(3),\GL_{2}(q))$}\label{Section:projective}

We consider $\GL_{2}(q)$ as a subgroup of $\SL_{3}(q)$ via the embedding $g\mapsto\diag(g,\det g^{-1})$.
We show that the pair $(\SL_{3}(q),\GL_{2}(q))$ does not admit any Gelfand triples. For each irreducible $\GL_{2}(q)$-representation $\pi$ we have to find an irreducible $\SL_{3}(q)$-representation that contains $\pi$ with multiplicity at least two.
We do not need to describe all the irreducible representations of $\SL_{3}(q)$ for this purpose. We do need the description of all the conjugacy classes of $\SL_{3}(q)$ and how the conjugacy classes of the embedded copy of $\GL_{2}(q)$ intersect with them.

Recall that $\rho$ and $\sigma$ are generators of $\bbF_{q}^{\times}$ and $\bbF_{q^{2}}^{\times}$ respectively. Let $\tau$ be the generator of $\bbF_{q^{3}}^{\times}$. 
We denote $t=q^2+q+1$, the number of elements in $\bbF^{\times}_{q^3}/\bbF^{\times}_{q}$.

The number of solutions of $X^3=1$ in $\bbF_{q}$ is denoted by $d$ and we have $d=\mathrm{gcd}(3,r)$. We denote $\omega=\rho^{r/d}$ so that $\omega^{d}=1$.
Furthermore, let $\theta\in\bbF_{q}^{\times}$ with $\theta^{3}\ne1$. 
In this section we also use a different way to describe the parameters of the irreducible representations. Indeed, in many occasions we work with subsets of $\bbZ$. For example, in this context the representations $V_{a}$ and $V_{a+r}$ are equivalent for $a\in\bbZ$. This is not a problem since we are not calculating the complete multiplicity table for the pair $(\SL_{3}(q),\GL_{2}(q))$.

Recall from \cite{SimpsonFrameSutherland} the conjugacy classes and the values of the characters as displayed in Tables \ref{SLcon} and 
\ref{SLchar}. Note that some of the representatives are in $\SL_{3}(q^2)$ or $\SL_{3}(q^3)$.
\begin{table}[h]
\scriptsize
\center
\begin{tabular}{|c c c|}
  \hline

    Conj.class & repr. & parameters \\ [0.5ex]
    \hline \hline
    $C_1(k)$ & $\begin{pmatrix}
    \omega^k & 0 & 0 &\\ 0 & \omega^k & 0 \\ 0 & 0 & \omega^k \end{pmatrix}$ & $1 \leq k \leq d$ \\
 & &  \\ 
$C_2(k)$ & $\begin{pmatrix} 
     \omega^k & 0 & 0 \\ 1 & \omega^k & 0 \\ 0 & 0 & \omega^k
\end{pmatrix}$ & $1 \leq k \leq d$  \\ 
& &  \\
$C_3(k,\ell)$ & $\begin{pmatrix}
     \omega^k & 0 & 0 \\ \theta^{\ell} & \omega^k & 0 \\ 0 & \theta^{\ell} & \omega^k
\end{pmatrix}$ & $1 \leq k,\ell \leq d$ \\
& & \\
$C_4(k)$ & $\begin{pmatrix}
     \rho^k & 0 & 0 \\ 0 & \rho^k & 0 \\ 0 & 0 & \rho^{-2k}
\end{pmatrix}$ & \makecell{$1 \leq k < r$\\ $k \neq 0\pmod{
\frac{r}{d}}$} \\
& & \\
$C_5(k)$ & $\begin{pmatrix}
     \rho^k & 0 & 0 \\ 1 & \rho^k & 0 \\ 0 & 0 & \rho^{-2k}
\end{pmatrix}$ & \makecell{$1 \leq k < r$ \\ $k \neq 0 \pmod{
\frac{r}{d}}$}\\
& & \\
$C_6(k,\ell)$ & $\begin{pmatrix}
     \rho^k & 0 & 0 \\ 0 & \rho^{\ell} & 0 \\ 0 & 0 & \rho^{-k-\ell}
\end{pmatrix}$ & \makecell{$1 \leq k < \ell \leq r$ \\ $k \neq -k-\ell \neq \ell$  $\pmod{r}$}\\
& &\\
$C_7(k)$ & $\begin{pmatrix}
     \rho^k & 0 & 0 \\ 0 & \sigma^{-k} & 0 \\ 0 & 0 & \sigma^{-qk}
\end{pmatrix}$ & \makecell{$1 \leq k < rs$ \\ $k \neq 0\pmod{s}$
\\ $C_7(k) = C_7(kq)$ $\pmod{rs}$} \\
& & \\
$C_8(k)$ & $\begin{pmatrix}
     \tau^k & 0 & 0 \\ 0 & \tau^{qk} & 0 \\ 0 & 0 & \tau^{kq^2}
\end{pmatrix}$ & \makecell{$1 \leq k < t$ \\ $k \neq 0\pmod{
\frac{t}{d}}$ \\ $C_8(k) = C_8(kq) = C_8(kq^2)$ $\pmod{t}$}\\
\hline
  \end{tabular}
  \caption{Conjugacy classes of $\SL_{3}(q)$}\label{SLcon}
\end{table}

\subsection{Restriction to $\GL_{2}(q)$}

Our aim is to restrict an irreducible representation $(\pi,V)$ of $\SL_{3}(q)$ to $\GL_{2}(q)$ and find an irreducible representation $(\tau,W)$ of $\GL_{2}(q)$ that occurs with multiplicity $>1$.
Let us denote 
$$[\pi|_{\GL_{2}(q)}:\tau]=\dim\Hom_{\GL_{2}(q)}(W,V).$$
To calculate this multiplicity we have to restrict characters of $\SL_{3}(q)$ to $\GL_{2}(q)$. For this we have to know in which conjugacy class of $\SL_{3}(q)$ a given conjugacy class of $\GL_{2}(q)$ sits. This depends on the conjugacy class and its parameter and the inclusions are displayed in Table \ref{embed}. 

\begin{table}[h]
    \centering
    \begin{tabular}{|c|c|}
    \hline & \\
         $c_1(k) \hookrightarrow \Bigg\{\makecell{C_1(k) ~\text{if}~ k = 0 \pmod{\frac{r}{d}} \\ C_4(k) ~\text{if}~ k \neq 0 \pmod{\frac{r}{d}}}$ & $\makecell{d \\ r-d}$ \\  & \\ \hline & \\ $c_2(k) \hookrightarrow \Bigg\{\makecell{C_2(k) ~\text{if}~ k = 0 \pmod{\frac{r}{d}} \\ C_5(k) ~\text{if}~ k \neq 0\pmod{\frac{r}{d}}}$ & $\makecell{d \\ r-d}$\\  &\\ \hline & \\ $ c_3(k,\ell) \hookrightarrow \left\{ \makecell{C_4(k) ~\text{if}~ 2k = -\ell \pmod{r} \\ C_4(\ell) ~\text{if}~ 2\ell = -k \pmod{r} \\ C_6(k,\ell) ~\text{if} ~ 2k \neq -\ell \pmod{r} \\ \text{and} ~ 2\ell \neq -k \pmod{r}} \right. $ & $\makecell{r-d \\ \\ \frac{r(r-1)}{2} - (r - d)}$ \\ & \\ \hline & \\ $c_4(k) \hookrightarrow C_7(-k)$ & $\frac{qr}{2}$ \\ &\\ \hline
         
    \end{tabular}
    \caption{Embedding of conjugacy classes of $\GL_{2}(q)$ into conjugacy classes of $\SL_{3}(q)$.}
    \label{embed}
\end{table}

The irreducible representations of $\SL_{3}(q)$ that we need are denoted by $\pi_{qs},\pi_{t^{(u)}}$ and $\pi_{rt^{(u)}}$ and their values on the conjugacy classes of $\SL_{3}(q)$ are displayed in Table \ref{SLchar}.

\begin{table}[h]
    \centering
    \begin{tabular}{c|c c c}
    Parameters & &$1 \leq u < r$ & $ \makecell{1 \leq u \leq rs \\ u \neq 0 \pmod{s} \\ u = qu \pmod{rs}}$ \\
      $\SL_{3}(q)$   & $\pi_{qs}$ & $\pi_{t^{(u)}}$ & $\pi_{rt^{(u)}}$ \\
      \hline \\
      $C_1(k)$ & $qs$ & $t\alpha_u(\omega^k)$ & $rt\alpha_u(\omega^k)$ \\
      $C_2(k)$ & $q$ & $s\alpha_u(\omega^k)$ & $-\alpha_u(\omega^k) $ \\
      $C_3(k,\ell)$ & $0$ & $\alpha_u(\omega^k)$ & $-\alpha_u(\omega^k) $ \\
      $C_4(k)$ & $s$ & $s\alpha_u(\rho^k)+\alpha_u(\rho^{-2k})$ & $r\alpha_u(\rho^k) $ \\
      $C_5(k)$ & $1$ & $\alpha_u(\rho^k)+\alpha_u(\rho^{-2k})$ & $-\alpha_u(\rho^k) $ \\
      $C_6(k,\ell)$ & $2$ & $\alpha_u(\rho^k)+\alpha_u(\rho^{\ell})+\alpha_u(\rho^{-k-\ell})$ & $0 $ \\
      $C_7(k)$ & $0$ & $\alpha_u(\rho^{k})$ & $-(\phi_u(\sigma^{-k})+\phi_u(\sigma^{-qk})) $ \\
      $C_8(k)$ & $-1$ & $0$ & $0 $ \\
      \\
      \hline
    \end{tabular}
    \caption{Partial character table for $\SL_{3}(q)$. Note that we have corrected the value on $C_{7}(k)$ in the last column.}
    \label{SLchar}
\end{table}

To calculate the multiplicities we refine \eqref{FormulaInnerProduct} according to the embeddings of the conjugacy classes of $\SL_{3}(q)$.

\begin{lemma}\label{SLformula} Let $(\pi,V)$ be an irreducible representation of $\SL_{3}(q)$ and $(\tau,W)$ and irreducible representation of $\GL_{2}(q)$. The multiplicity $[\pi|_{\GL_{2}(q)}:\tau]$ satisfies
\[
\begin{array}{c}
qsr^2 \left[\pi|_{\GL_{2}(q)}: \tau\right] =
\sum_{\substack{1 \leq k \leq r \\ k = 0\pmod{\frac{r}{d}}}} \chi_{\pi}(C_1(\tfrac{d}{r}k))\,\overline{\chi_{\tau}(c_1(k))} +
\sum_{\substack{1 \leq k \leq r \\ k \neq 0\pmod{\frac{r}{d}}}} \chi_{\pi}(C_4(k))\,\overline{\chi_{\tau}(c_1(k))}  \\[1ex]
+rs \left(
  \sum_{\substack{1 \leq k \leq r \\ k = 0\pmod{\frac{r}{d}}}} \chi_{\pi}(C_2(\tfrac{d}{r}k))\,\overline{\chi_{\tau}(c_2(k))} +
  \sum_{\substack{1 \leq k \leq r \\ k \neq 0\pmod{\frac{r}{d}}}} \chi_{\pi}(C_5(k))\,\overline{\chi_{\tau}(c_2(k))}
\right)  \\[1ex]
+qs \left(
  \sum_{\substack{1 \leq k \leq r \\ k \neq 0\pmod{\frac{r}{d}}}} \chi_{\pi}(C_4(k))\,\overline{\chi_{\tau}(c_3(k,-2k))} +
  \sum_{\substack{1 \leq k < \ell \leq r \\ \ell \neq -2k\pmod{r} \\ k \neq -2\ell\pmod{r}}} 
    \chi_{\pi}(C_6(k,\ell,-k-\ell))\,\overline{\chi_{\tau}(c_3(k,\ell))}
\right)  \\[1ex]
+ qr \left(
  \sum_{k \in \mathcal{X}} \chi_{\pi}(C_7(-k))\,\overline{\chi_{\tau}(c_4(k))}
\right).
\end{array}
\]
\end{lemma}

\begin{proof}
    This follows directly from \eqref{FormulaInnerProduct}, \eqref{FormulaMultiplicity} and embeddings of the conjugacy classes in Table \ref{embed}.
\end{proof}

\begin{lemma}
    \begin{itemize}
    \item $U_{0}$ occurs with multiplicity two in $\pi_{qs}$
    \item $U_{a}$ occurs in $\pi_{t^{(-a)}}$ with multiplicity two. 
    \item $V_{a}$ occurs in $\pi_{rt^{-a}}$ with multiplicity $d+1$.
    \item $W_{[a,b]}$ occurs $\pi_{rt^{(b-2a+cr)}}$ for $c$ with $b-2a+cr\ne0\pmod{s}$. The multiplicity is $1+d$ and even $2+d$ if moreover $b-a=0\pmod{r/d}$.    
    \item $X_{[n]}$ occurs with multiplicity two in $\pi_{rt^{(n)}}$.
\end{itemize}
\end{lemma}

\begin{proof}
    This follows from the calculations in Subsection \ref{SubSectionCalculationsSL}.
\end{proof}

\begin{corollary}
    The pair $(\SL_{3}(q),\GL_{2}(q))$ is not a Gelfand pair and it does not admit any Gelfand triples.
\end{corollary}

\subsection{Calculations}\label{SubSectionCalculationsSL}
Before we proceed with the calculations we note that the identity 
\[
\sum_{\substack{1 \leq j < k \leq r \\ j \neq -2k \pmod{r} \\ k \neq -2j \pmod{r}}} \alpha_a(\rho^{j+k})
= \frac{1}{2} \left( r^2 \, \mathbbm{1}_{a = 0 \pmod{r}} - r \, \mathbbm{1}_{2a = 0 \pmod{r}} \right)
- r \, \mathbbm{1}_{-a = 0 \pmod{r}}
+ d \, \mathbbm{1}_{-a = 0 \pmod{d}}.
\]
holds. Its proof uses the identities in Remark \ref{RemarkCalc} and counting over $\bbZ_{rs}\times\bbZ_{rs}$ first after which we subtract the terms that are not in the sum. This technique is used throughout our calculations below.

\subsubsection{Higer multiplicities in $\Ind U_a$}
First we take $a=0$, which makes $U_{0}$ the trivial representation, and we take $\pi_{qs}$,
\[
\begin{array}{c}
qsr^2[\pi_{qs}|_{\mathrm{GL}_2(q)} : U_0] = \\[1ex]
\sum_{\substack{1 \leq j \leq r \\ j = 0\pmod{\frac{r}{d}}}} qs +
\sum_{\substack{1 \leq j \leq r \\ j \neq 0\pmod{\frac{r}{d}}}} s + \\[1ex]
rs\left(
  \sum_{\substack{1 \leq j \leq r \\ j = 0\pmod{\frac{r}{d}}}} q +
  \sum_{\substack{1 \leq j \leq r \\ j \neq 0\pmod{\frac{r}{d}}}} 1
\right) + \\[1ex]
qs\left(
  \sum_{\substack{1 \leq j \leq r \\ j \neq 0\pmod{\frac{r}{d}}}} s +
  \sum_{\substack{1 \leq j < k \leq r \\ k \neq -2j\pmod{r} \\ j \neq -2k\pmod{r}}} 2
\right) + \\[1ex]
qr\left( \sum_{j \in \mathcal{X}} 0 \right) = \\[1ex]
dqs + s(r-d) + drsq + rs(r-d) + qs^2(r-d) + qs(r(r-1) - 2(r-d)) = 2qsr^2.
\end{array}
\]
Hence $[\pi_{qs}|_{\mathrm{GL}_2(q)} : U_0]=2$. For $a\ne0$ we take $\pi_{t^{(-a)}}$,
\[
\begin{array}{c}
qsr^2\left[\pi_{t^{-a}}|_{\mathrm{GL}_2(q)} : U_a\right] = \\[1ex]
\sum_{\substack{1 \leq j \leq r \\ j = 0\pmod{\frac{r}{d}}}}
  t\,\alpha_{-a}(\rho^j)\,\overline{\alpha_a(\rho^{2j})}
+ 
\sum_{\substack{1 \leq j \leq r \\ j \neq 0\pmod{\frac{r}{d}}}}
  \left(s\,\alpha_{-a}(\rho^j) + \alpha_{-a}(\rho^{-2j})\right)\,\overline{\alpha_a(\rho^{2j})} + \\[1ex]

rs \left(
  \sum_{\substack{1 \leq j \leq r \\ j = 0\pmod{\frac{r}{d}}}}
    s\,\alpha_{-a}(\rho^j)\,\overline{\alpha_a(\rho^{2j})} +
  \sum_{\substack{1 \leq j \leq r \\ j \neq 0\pmod{\frac{r}{d}}}}
    \left(\alpha_{-a}(\rho^j) + \alpha_{-a}(\rho^{-2j})\right)\,\overline{\alpha_a(\rho^{2j})}
\right) + \\[1ex]

qs \left(
  \sum_{\substack{1 \leq j \leq r \\ j \neq 0\pmod{\frac{r}{d}}}}
    \left(s\,\alpha_{-a}(\rho^j) + \alpha_{-a}(\rho^{-2j})\right)\,\overline{\alpha_a(\rho^{-j})}\right. \\
    +
\left.  \sum_{\substack{1 \leq j < k \leq r \\ k \neq -2j\pmod{r} \\ j \neq -2k\pmod{r}}}
    \left(\alpha_{-a}(\rho^j) + \alpha_{-a}(\rho^k) + \alpha_{-a}(\rho^{-j-k})\right)\,\overline{\alpha_a(\rho^{j+k})}
\right) + \\[1ex]

qr \left(
  \sum_{j \in \mathcal{X}} \alpha_{-a}(\rho^{-j})\,\overline{\alpha_a(\rho^j)}
\right).
\end{array}
\]
which we can rewrite as
\[
\begin{array}{c}
t\displaystyle\sum_{j=1}^{d} \alpha_{-3a\frac{r}{d}}(\rho^j)
+ \displaystyle\sum_{j=1}^r \left( s\,\alpha_{-3a}(\rho^j) + \alpha_0(\rho^j) \right)
- \left( \displaystyle\sum_{j=1}^{d} s\,\alpha_{-3a\frac{r}{d}}(\rho^j) + \alpha_{0\frac{r}{d}}(\rho^j) \right) + \\[1ex]

rs\left(
  s \displaystyle\sum_{j=1}^{d} \alpha_{-3a\frac{r}{d}}(\rho^j)
  + \displaystyle\sum_{j=1}^r \left( \alpha_{-3a}(\rho^j) + \alpha_0(\rho^j) \right)
  - \left( \displaystyle\sum_{j=1}^{d} \alpha_{-3a\frac{r}{d}}(\rho^j) + \alpha_{0\frac{r}{d}}(\rho^j) \right)
\right) + \\[1ex]

qs\left(
  \displaystyle\sum_{j=1}^r \left( s\,\alpha_0(\rho^j) + \alpha_{3a}(\rho^j) \right)
  - \left( \displaystyle\sum_{j=1}^{d} s\,\alpha_{0\frac{r}{d}}(\rho^j) + \alpha_{3a\frac{r}{d}}(\rho^j) \right)
  + \right. \\[1ex]

\left. \frac{1}{2} \left(
  \displaystyle\sum_{1 \leq j,k \leq r} \alpha_1(\rho^{(-2a)j - ak}) + \alpha_1(\rho^{(-2a)k - aj}) + \alpha_0(\rho^{j+k})
  - \left( \displaystyle\sum_{j=1}^r 2\alpha_{-3a}(\rho^j) + \alpha_0(\rho^j) \right)
\right) \right. - \\[1ex]

\left. \left( \displaystyle\sum_{j=1}^r 2\alpha_0(\rho^j) + \alpha_{3a}(\rho^j) \right)
+ \left( \displaystyle\sum_{j=1}^{d} 2\alpha_{0\frac{r}{d}}(\rho^j) + \alpha_{3a\frac{r}{d}}(\rho^j) \right)
\right) + \\[1ex]

\frac{1}{2}qr\left(
  \displaystyle\sum_{j=1}^{rs} \alpha_0(\rho^j)
  - \displaystyle\sum_{j=1}^r \alpha_0(\rho^j)
\right) = \\[1ex]

td - ds - d + rs^2d - 2rsd - qs^2d - qsd \\+ 3qsd + r + r^2s
- \frac{1}{2}qsr - \frac{1}{2}qr^2 + qs^2r + \frac{1}{2}qsr^2
- 2qrs + \frac{1}{2}qr^2s = 2qsr^2,
\end{array}
\]
where we use that $a \neq 0\pmod{r}$ and therefore $\sum_{1 \leq j,k \leq r}\alpha_1(\rho^{(-2a)j-ak})+\alpha_1(\rho^{(-2a)k-aj}) = 0$ for any choice of such $a$. It follows that $[\pi_{t^{-a}}|_{\mathrm{GL}_2(q)} : U_a]=2$.

\subsubsection{Higher multiplicities in $\Ind V_a$}

We take $\pi_{rt^{[-a]}}$ for $a \neq 0\pmod{r}$ and $\pi_{rt^{[-r]}}$ for $a=0 \pmod{r}$,
\[
\begin{array}{c}
qsr^2[\pi_{rt^{[-a]}}|_{\GL_{2}(q)}:V_{a}] = \\[1em]
\sum_{\substack{1 \leq j \leq r \\ j = 0\pmod{\frac{r}{d}}}}
rt\alpha_{-a}(\rho^j)\overline{q\alpha_{a}(\rho^{2j})}
+ 
\sum_{\substack{1 \leq j \leq r \\ j \neq 0 \pmod{\frac{r}{d}}}}
r\alpha_{-a}(\rho^j)\overline{q\alpha_{a}(\rho^{2j})}
~+ \\[1em]
qs\left(
\sum_{\substack{1 \leq j \leq r \\ j \neq 0 \pmod{\frac{r}{d}}}}
r\alpha_{-a}(\rho^j)\overline{\alpha_{a}(\rho^{-j})}
\right)
+ \\[1em]
qr\left(
\sum_{j \in \mathcal{X}}
\left(\phi_{-a}(\sigma^{j})+\phi_{-a}(\sigma^{qj})\right)
\overline{\alpha_{a}(\sigma^{sj})}
\right),
\end{array}
\]
which we can write as
\[
\begin{array}{c}
qrt\displaystyle\sum_{j=1}^{d}\alpha_{-3a\frac{r}{d}}(\rho^{j}) 
+ qr\left(\sum_{j=1}^r\alpha_{-3a}(\rho^j) 
- \sum_{j=1}^{d}\alpha_{-3a\frac{r}{d}}(\rho^{j})\right)
+ \\[1em]
qs\left(r\sum_{j=1}^r\alpha_{0}(\rho^j) 
- r\sum_{j=1}^{d}\alpha_{0\frac{r}{d}}(\rho^j)\right)
+ \\[1em]
qr\left(\sum_{j=1}^{rs}\phi_{(s+1)a}(\sigma^j) 
- \left(\sum_{j=1}^r\alpha_{-3a}(\rho^j)\right)\right)
= \\[1em]
dqrt - dqr + qs(r^2 - dr) = (d+1)qsr^2,
\end{array}
\]
since $d = \{1,3\} $ and $(s+1)a \neq 0 ~(rs)$ for $a \neq 0 ~(s)$. 
We conclude that $[\pi_{rt^{[-a]}}|_{\GL_{2}(q)}:V_{a}]=d+1$.

\subsubsection{Higher multiplicities in $\Ind W_{[a,b]}$}

We take $\pi_{rt^{[b-2a+cr]}}$ for a 
constant $c$ such that $b-2a+cr \neq 0 \pmod{s}$,
\[
\begin{array}{c}
qsr^2[\pi_{rt^{[b-2a+cr]}}|_{\GL_{2}(q)}:W_{[a,b]}] = \\[1em]
\sum_{\substack{1 \leq j \leq r \\ j = 0\pmod{\frac{r}{d}}}}
rt\alpha_{b-2a+cr}(\rho^j)
\overline{s\alpha_{a}(\rho^j)\alpha_{b}(\rho^j)}
+ 
\sum_{\substack{1 \leq j \leq r \\ j \neq 0\pmod{\frac{r}{d}}}}
r\alpha_{b-2a+cr}(\rho^j)
\overline{s\alpha_{a}(\rho^j)\alpha_{b}(\rho^j)}
~ - \\[1em]
rs\left(
\sum_{\substack{1 \leq j \leq r \\ j = 0\pmod{ \frac{r}{d}}}}
\alpha_{b-2a+cr}(\rho^j)
\overline{\alpha_{a}(\rho^j)\alpha_{b}(\rho^j)}
+
\sum_{\substack{1 \leq j \leq r \\ j \neq 0\pmod{\frac{r}{d}}}}
\alpha_{b-2a+cr}(\rho^j)
\overline{\alpha_{a}(\rho^j)\alpha_{b}(\rho^j)}
\right)
+ \\[1em]
qs\left(
\sum_{\substack{1 \leq j \leq r \\ j \neq 0\pmod{\frac{r}{d}}}}
r\alpha_{b-2a+cr}(\rho^j)
\overline{
\alpha_{a}(\rho^j)\alpha_{b}(\rho^{-2j})
+
\alpha_{a}(\rho^{-2j})\alpha_{b}(\rho^j)
}
\right),
\end{array}
\]
which we can write as
\[
\begin{array}{c}
srt\displaystyle\sum_{j=1}^{d}\alpha_{(cr-3a)\frac{r}{d}}(\rho^j) 
+ sr\left(\sum_{j=1}^r\alpha_{cr-3a}(\rho^j) 
- \sum_{j=1}^{d}\alpha_{(cr-3a)\frac{r}{d}}(\rho^j)\right)
- \\[1em]
rs\left(\sum_{j=1}^r\alpha_{cr-3a}(\rho^j)\right)
+ \\[1em]
qs\left(
r\sum_{j=1}^r\alpha_{3b-3a+cr}(\rho^j)
+ \alpha_{cr}(\rho^{j})
- r\sum_{j=1}^d 2\alpha_{(3a+cr)\frac{r}{d}}(\rho^j)
\right)
= \\[1em]
dsrt - dsr + qsr^2 - 2dqsr + qsr^2\mathbbm{1}_{3(b-a) = 0\pmod{r}}
= (d+1)qsr^2 + qsr^2\mathbbm{1}_{3(b-a) = 0\pmod{r}},
\end{array}
\]
from which it follows that $[\pi_{rt^{[b-2a+cr]}}|_{\GL_{2}(q)}:W_{[a,b]}]=$ is equal to $1 + d$ for any choice of $c$ such that $b-2a +cr \neq 0\pmod{s}$, or $2+d$ if $b-a$ is a multiple of $\frac{r}{d}$.
 
\subsubsection{Higher multiplicities in $\Ind X_{[n]}$}
We take $\pi_{rt^{[n]}}$, 

\[
\begin{array}{c}
qsr^2[\pi_{rt^{[n]}}|_{\GL_{2}(q)}:X_{[n]}] = \\[1em]
\sum_{\substack{1 \leq j \leq r \\ j = 0\pmod{\frac{r}{d}}}}
rt\alpha_n(\rho^j)\overline{r\alpha_n(\rho^j)}
+ 
\sum_{\substack{1 \leq j \leq r \\ j \neq 0\pmod{\frac{r}{d}}}}
r\alpha_n(\rho^j)\overline{r\alpha_n(\rho^j)}
+ \\[1em]
rs\left(
\sum_{\substack{1 \leq j \leq r \\ j = 0 \pmod{\frac{r}{d}}}}
\alpha_n(\rho^j)\overline{\alpha_n(\rho^j)}
+ 
\sum_{\substack{1 \leq j \leq r \\ j \neq 0 \pmod{\frac{r}{d}}}}
\alpha_n(\rho^j)\overline{\alpha_n(\rho^j)}
\right)
+ \\[1em]
qr\left(
\sum_{j \in \mathcal{X}}
\left(\phi_n(\sigma^j) + \phi_n(\sigma^{qj})\right)
\overline{\left(\phi_n(\sigma^j) + \phi_n(\sigma^{qj})\right)}
\right),
\end{array}
\]
which we can write as
\[
\begin{array}{c}
r^2t\sum_{j=1}^{d}\alpha_{0\frac{r}{d}}(\rho^j)
+ r^2\left(
\sum_{j=1}^r\alpha_0(\rho^j)
- \sum_{j=1}^{d}\alpha_{0\frac{r}{d}}(\rho^j)
\right)
+ \\[1em]
rs\left(\sum_{j=1}^r\alpha_0(\rho^j)\right)
+ \\[1em]
qr\left(
\sum_{j=1}^{rs} \phi_0(\sigma^j) + \phi_{rn}(\sigma^j)
- 2\sum_{j=1}^r\alpha_0(\rho^j)
\right)
= \\[1em]
dr^2t + r^3 - dr^2 + r^2s + qr^2s - 2qr^2 = (d+1)qsr^2.
\end{array}
\]
We conclude that $[\pi_{rt^{[n]}}|_{\GL_{2}(q)}:X_{[n]}]=d+1$.


\newpage

\appendix

\section{Calculations for Section \ref{Section:product}}\label{Appendix:GxG}

\begin{lemma}\label{Lemma:VW}
    The decomposition of $V_{a}\otimes W_{[b,c]}$ is described by
    \begin{itemize}
        \item $[V_{a}\otimes W_{[b,c]}:U_{a'}]=0$,
        \item $[V_{a}\otimes W_{[b,c]}:V_{b'}]=\mathbbm{1}_{2b'=2a+b+c}$,
        \item $[V_{a}\otimes W_{[b,c]}:W_{[c',d']}]=\mathbbm{1}_{c'+d'=2a+b+c}+
        \mathbbm{1}_{c'=a+b, d'=a+c}+
        \mathbbm{1}_{c'=a+c, d'=a+b}$,
        \item $[V_{a}\otimes W_{[b,c]}:X_{[n']}]=\mathbbm{1}_{\overline{n'}=2a+b+c}$.
    \end{itemize}
\end{lemma}

\begin{proof}
    We have $[V_{a}\otimes W_{[b,c]},U_{a'}]=0$ by Lemma \ref{sym} and Example \ref{ExampleTensorU} and $[V_{a}\otimes W_{[b,c]}:V_{b'}]=\mathbbm{1}_{2a+b+c=2b'}$ by Lemma \ref{sym} and Lemma \ref{Lemma:VV}.
    We continue with $[V_{a}\otimes W_{[b,c]}:W_{[c',d']}]$,
    \begin{multline*}
    qsr^2 [V_a \otimes W_{[b,c]} : W_{[c',d']}] = 
    qs^2 \sum_{j=1}^r \alpha_{2a+b+c-c'-d'}(\rho^{j})
    + qs \sum_{[j,k]\in\calW}\left(
    \alpha_{a+b-d'}(\rho^{j})\alpha_{a+c-c'}(\rho^{k})\right.\\
    \quad \left.+
    \alpha_{a+b-c'}(\rho^{j})\alpha_{a+c-d'}(\rho^{k})+
    \alpha_{a+c-d'}(\rho^{j})\alpha_{a+b-c'}(\rho^{k})+
    \alpha_{a+c-c'}(\rho^{j})\alpha_{a+b-d'}(\rho^{k})
    \right)\\
    = qs^2r \mathbbm{1}_{2a + b + c = c' + d'} - 
    2qsr(\mathbbm{1}_{2a+b+=c'+d'})+2qsr^2(\mathbbm{1}_{a+c=d',a+b=c'}+\mathbbm{1}_{a+b=d',a+c=c'})
    \\
    = (qsr^2) \mathbbm{1}_{2a + b + c = c' + d' }
    + (qsr^2) \mathbbm{1}_{a + b = c' , a + c = d' } + (qsr^2) \mathbbm{1}_{a + c = c', a + b = d'}.
    \end{multline*}
    Finaly,
    $$
    \begin{aligned}
    qsr^2 [V_a \otimes W_{[b,c]} : X_{[n']}] =& 
    qsr \sum_{j=1}^r \alpha_a(\rho^{2j}) \alpha_b(\rho^j) \alpha_c(\rho^j) \phi_{-n'}(\rho^j) \\
    =& qsr \sum_{j=1}^r \alpha_{2a + b + c - n'}(\rho^j) \\
    =& (qsr^2) \mathbbm{1}_{2a + b + c =\overline{n'}}
    \end{aligned}
    $$
    and we are done.
\end{proof}

\begin{lemma}\label{Lemma:VX}
    The decomposition of $V_{a}\otimes X_{n}$ is described by
    \begin{itemize}
        \item $[V_{a}\otimes X_{[n]}:U_{a'}]=0$,
        \item $[V_{a}\otimes X_{[n]}:V_{b'}]=\mathbbm{1}_{2a+\overline{n}=2b'}$,
        \item $[V_{a}\otimes X_{[n]}:W_{[c',d']}]=\mathbbm{1}_{2a+\overline{n}=c'+d'}$,
        \item $[V_{a}\otimes X_{[n]}:X_{[n']}]=\mathbbm{1}_{2a+\overline{n}=\overline{n'}}-\mathbbm{1}_{sa+n=n'}-\mathbbm{1}_{sa+n=qn}$.
    \end{itemize}
\end{lemma}

\begin{proof}
    The first three items are obtained from the first six rows of Table \ref{table:tensorproductdecomp} using Lemma \ref{sym} and 
    $$
    \begin{aligned}
    qsr^2[V_a \otimes X_{n}: X_{n'}] = & qr^2\sum_{j=1}^r\alpha_{2a+\overline{n}-\overline{n'}}(\rho^{j})\\
    & - qr\sum_{[j]\in\calX}(
    \phi_{sa+n-n'}(\sigma^j)+
    \phi_{sa+n-n'}(\sigma^{qj})+
    \phi_{sa+n-qn'}(\sigma^j)+
    \phi_{sa+n-n'}(\sigma^{qj}))\\
    =&qr^3\mathbbm{1}_{2a+\overline{n}=\overline{n'}}+2qr^2\mathbbm{1}_{2a+\overline{n}=\overline{n'}}-qsr^2\mathbbm{1}_{sa+n=n'}-qsr^2\mathbbm{1}_{sa+n=qn'}\\
    =&qsr^2\mathbbm{1}_{2a+\overline{n}=\overline{n'}}-qsr^2\mathbbm{1}_{sa+n=n'}-qsr^2\mathbbm{1}_{sa+n=qn'}
    \end{aligned}
    $$ 
    completes the proof.
\end{proof}

\begin{lemma}\label{Lemma:WW}
    The decomposition of $W_{[a,b]}\otimes W_{[c,d]}$ is described 
    \begin{itemize}
        \item $[W_{[a,b]}\otimes W_{[c,d]}:U_{a'}]=\mathbbm{1}_{a+c=a',b+d=a'}+\mathbbm{1}_{a+d=a',b+c=a'}$
        \item $[W_{[a,b]}\otimes W_{[c,d]}:V_{b'}]=\mathbbm{1}_{a+b+c+d=2b'}+\mathbbm{1}_{a+c=b',b+d=b'}+\mathbbm{1}_{a+d=b',b+c=b'}$
        \item $[W_{[a,b]}\otimes W_{[c,d]}:W_{[c',d']}]=\mathbbm{1}_{a+b+c+d=c'+d'}
        + \mathbbm{1}_{a+c=c', b+d=d'}
        + \mathbbm{1}_{a+c=d', b+d=c'}
        + \mathbbm{1}_{a+d=c', b+c=d'}
        + \mathbbm{1}_{a+d=d', b+c=c'}$
        \item $[W_{[a,b]}\otimes W_{[c,d]}:X_{[n']}]=\mathbbm{1}_{a+b+c+d=\overline{n'}}$
    \end{itemize}
\end{lemma}

\begin{proof}
    The first two items are obtained from the first rows of Table \ref{table:tensorproductdecomp} using Lemma \ref{sym}. Next,
    \begin{multline*}
    qsr^2[W_{[a,b]} \otimes W_{[c,d]}: W_{[c',d']}] =  s^3\sum_{j=1}^r\alpha_{a+b+c+d-c'-d'}(\rho^j) +rs\sum_{j=1}^r\alpha_{a+b+c+d-c'-d'}(\rho^j)\\
    + qs\left(\sum_{j=1}^r\sum_{k=1}^r \alpha_{a+c-c'}(\rho^j)\alpha_{b+d-d'}(\rho^k) + \alpha_{a+c-d'}(\rho^j)\alpha_{b+d-c'}(\rho^k) + \right.\\
    + \alpha_{a+d-c'}(\rho^j)\alpha_{b+c-d'}(\rho^k) + \alpha_{a+d-d'}(\rho^j)\alpha_{b+c-c'}(\rho^k) - 4\left.\sum_{j=1}^r\alpha_{a+b+c+d-c'-d'}(\rho^j)\right)\\
    = (s^3r+r^2s-4qsr)\mathbbm{1}_{a+b+c+d-c'-d' = 0} + 
    (qsr^2)\mathbbm{1}_{a+c-c' = b+d-d' = 0} + (qsr^2)\mathbbm{1}_{a+c-d' = b+d-c' = 0}\\
    +(qsr^2)\mathbbm{1}_{a+d-c' = b+c-d' = 0} + (qsr^2)\mathbbm{1}_{a+d-d' = b+c-c' = 0 }\\
    =qsr^2\left(
    \mathbbm{1}_{a+b+c+d=c'+d'}
    + \mathbbm{1}_{a+c=c', b+d=d'}
    + \mathbbm{1}_{a+c=d', b+d=c'}
    + \mathbbm{1}_{a+d=c', b+c=d'}
    + \mathbbm{1}_{a+d=d', b+c=c'}
    \right),
    \end{multline*}
    and finally
    $$\begin{aligned}
    qsr^2[W_{[a,b]} \otimes W_{[c,d]}: X_{[n']}] =&  s^2r\sum_{j=1}^r\alpha_{a+b+c+d-n'}(\rho^j) -rs\sum_{j=1}^r\alpha_{a+b+c+d-n'}(\rho^j) \\
    =& qsr\sum_{j=1}^r\alpha_{a+b+c+d-n'}(\rho^j) = (qsr^2)\mathbbm{1}_{a+b+c+d-\overline{n'} = 0}
    \end{aligned}
    $$
    as desired.
\end{proof}

\begin{lemma}\label{Lemma:WX}
    The decomposition of $W_{[a,b]}\otimes X_{[n]}$ is described by:
    \begin{itemize}
        \item $[W_{[a,b]}\otimes X_{[n]}:U_{a'}]= 0$
        \item $[W_{[a,b]}\otimes X_{[n]}:V_{b'}]= \mathbbm{1}_{a+b+\overline{n}=2b'}$
        \item $[W_{[a,b]}\otimes X_{[n]}:W_{[c',d']}]= \mathbbm{1}_{a+b+\overline{n}=c'+d'}$
        \item $[W_{[a,b]}\otimes X_{[n]}:X_{[n']}]=\mathbbm{1}_{a+b+\overline{n}=\overline{n'}}$
    \end{itemize}
\end{lemma}

\begin{proof}
   The first three items are obtained from the first rows of Table \ref{table:tensorproductdecomp} using Lemma \ref{sym}, so
    $$\begin{aligned}
    qsr^2[W_{[a,b]} \otimes X_{[n]}: X_{[n']}] = & sr^2\sum_{j=1}^r\alpha_{a+b+n-n'}(\rho^j) +rs\sum_{j=1}^r\alpha_{a+b+n-n'}(\rho^j)\\
    = &  qsr\sum_{j=1}^r\alpha_{a+b+n-n'}(\rho^j) = (qsr^2)\mathbbm{1}_{a+b+\overline{n}=\overline{n'}}
    \end{aligned}
    $$
proves the claim.
\end{proof}

\begin{lemma}\label{Lemma:XX}
    The decomposition of $X_{[n]}\otimes X_{[m]}$ is described by:
    \begin{itemize}
        \item $[X_{[n]}\otimes X_{[m]}:U_{a'}]= \mathbbm{1}_{n+m=sa'}$
        \item $[X_{[n]}\otimes X_{[m]}:V_{b'}]= \mathbbm{1}_{\overline{n}+\overline{m}=2b'}-\mathbbm{1}_{n+m=sb}-\mathbbm{1}_{n+qm=sb}$
        \item $[X_{[n]}\otimes X_{[m]}:W_{[c',d']}]= \mathbbm{1}_{\overline{n}+\overline{m}=c'+d'}$
        \item $[X_{[n]}\otimes X_{[m]}:X_{[n']}]= \mathbbm{1}_{\overline{n}+\overline{m}=\overline{n'}}
        - \mathbbm{1}_{n+m=n'} - \mathbbm{1}_{qn+m=n'}
        - \mathbbm{1}_{n+m=qn'} - \mathbbm{1}_{qn+m=qn'}$.
    \end{itemize}
\end{lemma}

\begin{proof}
    The first three items are obtained from the first rows of Table \ref{table:tensorproductdecomp} using Lemma \ref{sym}, so
    $$\begin{aligned}
    qsr^2[X_{[n]} \otimes X_{[m]}: X_{[n']}] = &    
    \end{aligned}
    $$
$$r^3\sum_{j=1}^r\alpha_{n+m-n'}(\rho^j) - rs\sum_{j=1}^r\alpha_{n+m-n'}(\rho^j) -$$ $$qr\left(\sum_{j=1}^{rs}\phi_{n+m-n'}(\sigma^j) + \phi_{qn+m-n'}(\sigma^j) + \phi_{n+qm-n'}(\sigma^j) + \right. $$ $$  \left.\phi_{n+m-qn'}(\sigma^j) -4\sum_{j=1}^r\alpha_{n+m-n'}(\rho^j) \right) =$$
$$ (r^4-r^2s+4qr^2)\mathbbm{1}_{\overline{n}+\overline{m}-\overline{n'} = 0 } - (qsr^2)\mathbbm{1}_{n+m-n' = 0} - (qsr^2)\mathbbm{1}_{qn+m-n' = 0} - $$ $$ (qsr^2)\mathbbm{1}_{n+qm-n' = 0 } - (qsr^2)\mathbbm{1}_{n+m-qn' = 0 } =$$
$$ (qsr^2)\mathbbm{1}_{\overline{n}+\overline{m}-\overline{n'} = 0 } - (qsr^2)\mathbbm{1}_{n+m-n' = 0} - (qsr^2)\mathbbm{1}_{qn+m-n' = 0 } - $$ $$ (qsr^2)\mathbbm{1}_{n+qm-n' = 0 } - (qsr^2)\mathbbm{1}_{n+m-qn' = 0 },$$
since $r^4-r^2s+4qr^2 = (q-3)qr^2+4qr^2 = qsr^2$.
Thus $[X_{[n]} \otimes X_{[m]}: X_{[n']}] = 1$ if and only if $\overline{n}+\overline{m}-\overline{n'} = 0$ $\pmod{r}$ and none of the following hold: $n+m-n' = 0$ $\pmod{rs}$, $qn+m-n' = 0$ $\pmod{rs}$, $n+qm-n' = 0$ $\pmod{rs}$ or $n+m-qn' = 0$ $\pmod{rs}$. Note again that only one of these last four conditions can apply at the same time.
\end{proof}

\newpage

\section{Examples for Section \ref{Section:product}}

In this appendix we illustrate the results of Section \ref{Section:product} by giving some explicit decompositions of restrictions and inductions from $\GL_{2}(q)\times\GL_{2}(q)$ to $\diag(\GL_{2}(q))$ and back. In this appendix we work with parameters in subsets of $\bbZ$.

\begin{example}
    Note that $\dim V_{a}\otimes V_{b}=q^2$. Let us verify that the dimensions of the irreducible submodules add up to $q^2$. We have to make a case distinction between $q$ odd and even.

    Suppose that $q$ is odd. The only one dimensional irreducible submodule is $U_{a+b}$. The only irreducible submodules of dimension $q$ are $V_{a+b}$ and $V_{a+b+r/2}$. There are $(r-2)/2$ irreducible submodules of dimension $q+1$,
    $$W_{[a+b+1,a+b-1]},\quad W_{[a+b+2,a+b-2]} ,\quad\ldots, \quad W_{[a+b+(\frac{r}{2}-1),a+b -(\frac{r}{2}-1)]},$$
    and $(s-2)/2$ submodules of dimension $q-1$,
    $$X_{[2(a+b)+(a+b+1)r]},\quad X_{[2(a+b)+(a+b+2)r]},\quad\ldots,\quad X_{[2(a+b)+(a+b+\frac{1}{2}s-1)r]}$$
    and the dimension count is $1 + 2q + \frac{r-2}{2}(q+1) + \frac{s-2}{2}(q-1) = q^2$
    as expected.

    If $q$ is even, then the only irreducible submodules of dimension $1$ and $q$ are $U_{a+b}$ and $V_{a+b}$. There are $(r-1)/2$ irreducible submodules of dimension $q+1$, namely
    $$W_{[a+b+1,a+b-1]},\quad W_{[a+b+2,a+b-2]} ,\quad\ldots, \quad W_{[a+b+\frac{r-1}{2},a+b -\frac{r-1}{2}]},$$
    and $(s-1)/2$ irreducible submodules of dimension $q-1$,
    $$X_{[2(a+b)+(a+b+1)r]},\quad X_{[2(a+b)+(a+b+2)r]},\quad\ldots,\quad X_{[2(a+b)+(a+b+\frac{s-1}{2} )r]}$$
    which gives the dimension count $1 + q + \frac{r-1}{2}(q+1) + \frac{s-1}{2}(q-1) = q^{2}$ 
    as expected.
\end{example}
    
\begin{example}
    The irreducible $G$-module $U_{a'}$ occurs in $r$ the irreducible $G\times G$-modules $V_{a}\otimes V_{b}$ for $a+b=a'\pmod{r}$, 
    $$V_{a'}\otimes V_{r},\quad V_{a'+1}\otimes V_{r-1},\ldots,V_{a'+r-1}\otimes V_{1}.$$
    If $q$ is odd and $n'$ is even, the irreducible $G$-module $X_{[n']}$ occurs in the $2r$ modules $V_{a}\otimes V_{b}$ with $2a+2b=n'\pmod{r}$,
    $$V_{0}\otimes V_{\frac{1}{2}n'},\quad V_{1}\otimes V_{\frac{1}{2}n'-1},\ldots,V_{r-1}\otimes V_{\frac{1}{2}n'-(r-1)}$$ and $$V_{0}\otimes V_{\frac{1}{2}n'+\frac{1}{2}r},\quad V_{1}\otimes V_{\frac{1}{2}n'-1+\frac{1}{2}r},\ldots,V_{r-1}\otimes V_{\frac{1}{2}n'-(r-1)+\frac{1}{2}r}.$$ If $n'$ is odd, the irreducible $G$-module $X_{n'}$ does not occur in any module of the type $V_a \otimes V_b$.
    If $q$ is even,then  $2$ has a multiplicative inverse modulo $r$ and the irreducible $G$-module $X_{[n']}$ occurs in the $r$ modules $V_{a}\otimes V_{b}$ with $2a+2b=n'\pmod{r}$,
    $$V_{0}\otimes V_{\frac{r+1}{2}n'},\quad V_{1}\otimes V_{\frac{r+1}{2}n'-1},\ldots,V_{r-1}\otimes V_{\frac{r+1}{2}n'-(r-1)}.$$
\end{example}

\begin{example}
    Note that for each $V_{a}\otimes W_{[b,c]}$ there is a module $W_{[c',d']}$ for which the multiplicity is two. Indeed, take $W_{[a+b,a+c]}$.
    Conversely, each module $W_{[c',d']}$ occurs in the decomposition of some tensor product $V_{a}\otimes W_{[b,c]}$ with multiplicity two, namely $$V_0 \otimes W_{[c',d']},\quad V_1 \otimes W_{[c'-1,d'-1]},\quad\ldots,\quad V_{r-1}\otimes W_{[c'(r-1),d'-(r-1)]}.$$

    Similarly, if $a+c = b+d$ and $q$ odd, then $W_{[a,b]}\otimes W_{[c,d]}$ contains $V_{a+c}$ with multiplicity two and $V_{b'}$ occurs in $$W_{[b',b'+1]}\otimes W_{[0,-1]},\quad W_{[b',b'+2]}\otimes W_{[0,-2]},\quad\ldots,\quad W_{[b',b+(\frac{1}{2}r-2)]}\otimes W_{[0,-(\frac{1}{2}r-2)]} \quad W_{[b',b'+(\frac{1}{2}r-1)]}\otimes W_{[0,-(\frac{1}{2}r-1)]}$$ 
    $$W_{[b'+1,b'+2]}\otimes W_{[-1,-2]},\quad W_{[b'+1,b'+3]}\otimes W_{[-1,-3]},\quad\ldots,\quad W_{[b'+1,b+(\frac{1}{2}r-1)]}\otimes W_{[-1,-(\frac{1}{2}r-1)]} $$ 
    $$\ldots$$
    $$W_{[b'+(\frac{1}{2}r-3),b'+(\frac{1}{2}r-2)]} \otimes W_{[-(\frac{1}{2}r-3),-(\frac{1}{2}r-2)]}, \quad W_{[b'+(\frac{1}{2}r-3),b'+(\frac{1}{2}r-1)]} \otimes W_{[-(\frac{1}{2}r-3),-(\frac{1}{2}r-1)]}$$
    $$W_{[b'+(\frac{1}{2}r-2),b'+(\frac{1}{2}r-1)]} \otimes W_{[-(\frac{1}{2}r-2),-(\frac{1}{2}r-1)]}$$
    with multiplicity two. 
\end{example}

\begin{example} Here we give the explicit decomposition of $V_{a}\otimes W_{[b,c]}$. 
If $q$ is odd and $b+c$ is even, then
    \begin{itemize}
        \item $2$ $q$-dimensional irreducible submodules 
        $$V_{a+\frac{1}{2}(b+c)},\quad V_{a+\frac{1}{2}(b+c)+\frac{1}{2}r}$$
        \item $\frac{r-2}{2}+1$ $(q+1)$-dimensional irreducible submodules 
        $$W_{[a+\frac{1}{2}(b+c) +1,a+\frac{1}{2}(b+c)-1]},\quad W_{[a+\frac{1}{2}(b+c) +2,a+\frac{1}{2}(b+c)-2]}, \quad\ldots, \quad W_{[a+\frac{1}{2}(b+c) +\frac{1}{2}r-1,a+\frac{1}{2}(b+c)-(\frac{1}{2}r-1)]},$$  with $W_{[a+b,a+c]}$ occurring twice.
        \item $\frac{s-2}{2}$ $(q-1)$-dimensional irreducible submodules 
        $$X_{[2a+b+c+(a+\frac{1}{2}(b+c)+1)r]},\quad X_{[2a+b+c+(a+\frac{1}{2}(b+c)+2)r]},\quad\ldots,\quad X_{[2a+b+c+(a+\frac{1}{2}(b+c)+\frac{s-2}{2})r]}$$
    \end{itemize}
    and the dimension count reads $2q+\frac{r}{2}(q+1) + \frac{r}{2}(q-1)=q(q+1)$ as expected.\\
    If $q$ is odd and $b+c$ is odd as well, then the module $V_a \otimes W_{[b,c]}$ decomposes as a direct sum of the
    \begin{itemize}

        \item $\frac{r}{2}+1$ $(q+1)$-dimensional irreducible submodules 
        $$W_{a+\frac{1}{2}(b+c+1) ,a+\frac{1}{2}(b+c-1)},\quad W_{[a+\frac{1}{2}(b+c+1) +1,a+\frac{1}{2}(b+c+1)-1]}, \quad\ldots, \quad W_{[a+\frac{1}{2}(b+c+1) +r-1,a+\frac{1}{2}(b+c-1)-(r-1)]},$$ with $W_{[a+b,a+c]}$ occurring twice.
        \item $\frac{s}{2}$ $(q-1)$-dimensional irreducible submodules 
        $$X_{[2a+b+c+(a+\frac{1}{2}(b+c+1))r]},\quad X_{[2a+b+c+(a+\frac{1}{2}(b+c+1)+1)r]},\quad\ldots,\quad X_{[2a+b+c+(a+\frac{1}{2}(b+c+1)+\frac{s-2}{2})r]}$$
    \end{itemize}
    and the dimension count reads $\frac{s}{2}(q+1) +\frac{s}{2}(q-1) = q(q+1)$.\\

    If $q$ is even, then the module $V_a \otimes W_{[b,c]}$ decomposes as a direct sum of the
    \begin{itemize}
       
         \item   $q$-dimensional irreducible submodule 
        $$V_{a+\frac{r+1}{2}(b+c)}$$
        \item $\frac{r-1}{2}+1$ $(q+1)$-dimensional irreducible submodules 
        $$W_{[a+\frac{r+1}{2}(b+c)+1 ,a+\frac{r+1}{2}(b+c)-1]},\quad W_{[a+\frac{r+1}{2}(b+c) +2,a+\frac{r+1}{2}(b+c)-2]}, \quad\ldots, \quad W_{[a+\frac{r+1}{2}(b+c) +\frac{r-1}{2},a+\frac{r+1}{2}(b+c)-(\frac{r-1}{2})]},$$ with $W_{[a+b,a+c]}$ occurring twice.
        \item $\frac{s-1}{2}$ $(q-1)$-dimensional irreducible submodules 
        $$X_{[2a+b+c+(a+\frac{rs+1}{2}(b+c)+1)r]},\quad X_{[2a+b+c+(a+\frac{rs+1}{2}(b+c)+2)r]},\quad\ldots,\quad X_{[2a+b+c+(a+\frac{rs+1}{2}(b+c)+\frac{s-1}{2})r]}$$
    \end{itemize}
    and the dimension count reads $q+\frac{q}{2}(q+1) +\frac{q}{2}(q-1) = q(q+1)$.\\
\end{example}

\begin{example}
    
    The irreducible $G$-module $U_{a'}$ does not occur in any of the irreducible $G\times G$-modules $V_{a}\otimes W_{[b,c]}$.
    If $q$ is odd and $n'$ is even, the irreducible $G$-module $X_{[n']}$ occurs in the $\frac{1}{2}r(r-2)$ modules $V_{a}\otimes W_{[b,c]}$ and the $\frac{1}{2}r(r-2)$ modules $W_{[b,c]} \otimes V_a$ with $2a+b+c=n'$,
    $$V_{0}\otimes W_{[\frac{1}{2}n'+1,\frac{1}{2}n'-1]},\quad V_{0}\otimes W_{[\frac{1}{2}n'+2,\frac{1}{2}n'-2]},\quad\ldots, \quad V_{0}\otimes W_{[\frac{1}{2}n'+\frac{1}{2}r-1,\frac{1}{2}n'-(\frac{1}{2}r-1)]}$$
    $$V_{1}\otimes W_{[\frac{1}{2}n'-2+1,\frac{1}{2}n'-2-1]},\quad V_{1}\otimes W_{[\frac{1}{2}n'-2+2,\frac{1}{2}n'-2-2]},\quad\ldots, \quad V_{1}\otimes W_{[\frac{1}{2}n'-2+\frac{1}{2}r-1,\frac{1}{2}n'-2-(\frac{1}{2}r-1)]}$$
    $$\ldots$$
    $$V_{r-1}\otimes W_{[\frac{1}{2}n'-2(r-1)+1,\frac{1}{2}n'-2(r-1)-1]},\quad V_{r-1}\otimes W_{[\frac{1}{2}n'-2(r-1)+2,\frac{1}{2}n'-2(r-1)-2]},\quad\ldots, $$$$ V_{r-1}\otimes W_{[\frac{1}{2}n'-2(r-1)+\frac{1}{2}r-1,\frac{1}{2}n'-2(r-1)-(\frac{1}{2}r-1)]},$$ and you get the $W_{[b,c]} \otimes V_a$ modules by "swapping" the ones given above. \\

    If $q$ is odd and $n'$ is odd, the irreducible $G$-module $X_{[n']}$ occurs in the $\frac{1}{2}r^2$ modules $V_{a}\otimes W_{[b,c]}$ and the $\frac{1}{2}r^2$ modules $W_{[b,c]} \otimes V_a$ with $2a+b+c=n'$,
    $$V_{0}\otimes W_{[\frac{1}{2}(n'+1),\frac{1}{2}(n'-1)]},\quad V_{0}\otimes W_{[\frac{1}{2}(n'+1)+1,\frac{1}{2}(n'-1)-1]},\quad\ldots, \quad V_{0}\otimes W_{[\frac{1}{2}(n'+1)+\frac{1}{2}r-1,\frac{1}{2}(n'-1)-(\frac{1}{2}r-1)]}$$
    $$V_{1}\otimes W_{[\frac{1}{2}(n'+1)-2,\frac{1}{2}(n'-1)-2]},\quad V_{1}\otimes W_{[\frac{1}{2}(n'+1)-2+1,\frac{1}{2}(n'-1)-2-1]},\quad\ldots, \quad V_{1}\otimes W_{[\frac{1}{2}(n'+1)-2+\frac{1}{2}r-1,\frac{1}{2}(n'-1)-2-(\frac{1}{2}r-1)]}$$
    $$\ldots$$
    $$V_{r-1}\otimes W_{[\frac{1}{2}(n'+1)-2(r-1),\frac{1}{2}(n'-1)-2(r-1)]},\quad V_{r-1}\otimes W_{[\frac{1}{2}(n'+1)-2(r-1)+1,\frac{1}{2}(n'-1)-2(r-1)-1]},\quad\ldots, $$$$ V_{r-1}\otimes W_{[\frac{1}{2}(n'+1)-2(r-1)+\frac{1}{2}r-1,\frac{1}{2}(n'-1)-2(r-1)-(\frac{1}{2}r-1)]},$$ and you get the $W_{[b,c]} \otimes V_a$ modules by "swapping" the ones given above. \\

    If $q$ is even, the irreducible $G$-module $X_{[n']}$ occurs in the $\frac{1}{2}r(r-1)$ modules $V_{a}\otimes W_{[b,c]}$ and the $\frac{1}{2}r(r-1)$ modules $W_{[b,c]} \otimes V_a$ with $2a+b+c=n'$,
    $$V_{0}\otimes W_{[\frac{r+1}{2}n'+1,\frac{r+1}{2}n'-1]},\quad V_{0}\otimes W_{[\frac{r+1}{2}n'+2,\frac{r+1}{2}n'-2]},\quad\ldots, \quad V_{0}\otimes W_{[\frac{r+1}{2}n'+\frac{r-1}{2},\frac{r+1}{2}n'-(\frac{r-1}{2})]}$$
    $$V_{1}\otimes W_{[\frac{r+1}{2}n'-2+1,\frac{r+1}{2}n'-2-1]},\quad V_{1}\otimes W_{[\frac{r+1}{2}n'-2+2,\frac{r+1}{2}n'-2-2]},\quad\ldots, \quad V_{1}\otimes W_{[\frac{r+1}{2}n'-2+\frac{r-1}{2},\frac{r+1}{2}n'-2-(\frac{r-1}{2})]}$$
    $$\ldots$$
    $$V_{r-1}\otimes W_{[\frac{r+1}{2}n'-2(r-1)+1,\frac{r+1}{2}n'-2(r-1)-1]},\quad V_{r-1}\otimes W_{[\frac{r+1}{2}n'-2(r-1)+2,\frac{r+1}{2}n'-2(r-1)-2]},\quad\ldots, $$$$ V_{r-1}\otimes W_{[\frac{r+1}{2}n'-2(r-1)+\frac{r-1}{2},\frac{r+1}{2}n'-2(r-1)-(\frac{r-1}{2})]},$$ and you get the $W_{[b,c]} \otimes V_a$ modules by "swapping" the ones given above.

\end{example}

\begin{example}
    There are no irreducible $G\times G$-modules $V_{a}\otimes X_{[n]}$ that contain $U_{a'}$.\\
    If $q$ is odd and $n'$ is even, then there are $\frac{1}{2}r(s-4)$ irreducible $G\times G$-modules $V_{a}\otimes X_{[n]}$ and $\frac{1}{2}r(s-4)$ irreducible $G\times G$-modules $X_{[n]}\otimes V_{a}$ that contain $X_{[n']}$,
    $$V_a\otimes X_{[n'+(\frac{1}{2}n'-a+1)r]}, \quad V_a\otimes X_{[n'+(\frac{1}{2}n'-a+2)r]}, \quad \ldots, \quad V_a \otimes X_{[n'+(\frac{1}{2}n'-a+\frac{s-2}{2})r]}$$
    for any $1 \leq a \leq r$, where we omit $V_a \otimes X_{[n'-sa]}$ for each $a$. Again, we swap to get the irreps of type $X_{[n]} \otimes V_a$. \\
    If $q$ is odd and $n'$ is odd as well, then there are $\frac{1}{2}r(s-2)$ irreducible $G\times G$-modules $V_{a}\otimes X_{[n]}$ and $\frac{1}{2}r(s-2)$ irreducible $G\times G$-modules $X_{[n]}\otimes V_{a}$ that contain $X_{[n']}$,
    $$V_a\otimes X_{[n'+(\frac{1}{2}(n'+1)-a)r]}, \quad V_a\otimes X_{[n'+(\frac{1}{2}(n'+1)-a+1)r]}, \quad \ldots, \quad V_a \otimes X_{[n'+(\frac{1}{2}(n'+1)-a+\frac{s-2}{2})r]}$$
    for any $1 \leq a \leq r$, where we omit $V_a \otimes X_{[n'-sa]}$ for each $a$. Again, we swap to get the irreps of type $X_{[n]} \otimes V_a$. \\
    If $q$ is even, then there are $\frac{1}{2}r(s-3)$ irreducible $G\times G$-modules $V_{a}\otimes X_{[n]}$ and $\frac{1}{2}r(s-3)$ irreducible $G\times G$-modules $X_{[n]}\otimes V_{a}$ that contain $X_{[n']}$,
    $$V_a\otimes X_{[n'+(\frac{rs+1}{2}(n')-a+1)r]}, \quad V_a\otimes X_{[n'+(\frac{rs+1}{2}(n')-a+2)r]}, \quad \ldots, \quad V_a \otimes X_{[n'+(\frac{rs+1}{2}(n')-a+\frac{s-1}{2})r]}$$
    for any $1 \leq a \leq r$, where we omit $V_a \otimes X_{[n'-sa]}$ for each $a$. Again, we swap to get the irreps of type $X_{[n]} \otimes V_a$.
\end{example}

In the following we list certain tensor products that contain a given representation $U_{a'}$. Given a representation $X_{[n']}$ we do not list all the tensor products of these types that contain $X_{[n']}$ upon restriction, since this is complicated and needs more notation. Instead we count the number of these tensor porduct representations, see also Examples \ref{exampleWW}, \ref{exampleWX} and \ref{exampleXX}.

\begin{example}
    There are $\frac{1}{2}r(r-1)$ irreducible $G\times G$-modules $W_{[a,b]}\otimes W_{[c,d]}$ that contain $U_{a'}$,
    $$W_{[0,1]} \otimes W_{[a',a'-1]}, \quad W_{[0,2]} \otimes W_{[a' ,a' -2]}, \quad \ldots, \quad W_{[0,r-2]} \otimes W_{[a',a'-(r-2)]},\quad W_{[0,r-1]} \otimes W_{[a',a'-(r-1)]}$$
    $$W_{[1,2]} \otimes W_{[a'-1,a'-2]}, \quad W_{[1,3]} \otimes W_{[a'-1 ,a' -3]}, \quad \ldots, \quad W_{[0,r-1]} \otimes W_{[a',a'-(r-1)]}$$
    $$\ldots$$
    $$ W_{[r-3,r-2]} \otimes W_{[a'-(r-3),a'-(r-2)]} \quad W_{[r-3,r-1]} \otimes W_{[a'-(r-3),a'-(r-1)]}$$
    $$W_{[r-2,r-1]} \otimes W_{[a'-(r-2),a'-(r-1)]}.$$ 
    If $q$ is odd and $n'$ is even, there are $\frac{1}{4}(r^3-2r^2+2r)$ irreducible $G\times G$-modules $W_{[a,b]}\otimes W_{[c,d]}$ that contain $X_{[n']}$.
    If $q$ is odd and $n'$ is odd as well, there are $\frac{1}{4}(r^3-2r^2)$ irreducible $G\times G$-modules $W_{[a,b]}\otimes W_{[c,d]}$ that contain $X_{[n']}$.
    If $q$ is even, there are $\frac{1}{4}(r^3-2r^2+r)$ irreducible $G\times G$-modules $W_{[a,b]}\otimes W_{[c,d]}$ that contain $X_{[n']}$.
    
\end{example}

\begin{example}
    There are no irreducible $G\times G$-modules $W_{[a,b]}\otimes X_{[n]}$ that contain $U_{a'}$. \\
    
    If $q$ is odd and $n'$ is even, there are $\frac{1}{4}(r^2s-rs-r^2+2r)$ irreducible $G\times G$-modules $W_{[a,b]}\otimes X_{[n]}$ and $\frac{1}{4}(r^2s-rs-r^2+2r)$ irreducible $G\times G$-modules $X_{[n]}\otimes W_{[a,b]}$ that contain $X_{[n']}$.
    If $q$ is odd and $n'$ is odd as well, there are $\frac{1}{4}(r^2s-rs-r^2)$ irreducible $G\times G$-modules $W_{[a,b]}\otimes X_{[n]}$ and $\frac{1}{4}(r^2s-rs-r^2)$ irreducible $G\times G$-modules $X_{[n]}\otimes W_{[a,b]}$ that contain $X_{[n']}$.
    If $q$ is even, there are $\frac{1}{4}(r^2s-rs-r^2+r)$ irreducible $G\times G$-modules $W_{[a,b]}\otimes X_{[n]}$ and $\frac{1}{4}(r^2s-rs-r^2+r)$ irreducible $G\times G$-modules $X_{[n]}\otimes W_{[a,b]}$ that contain $X_{[n']}$.
\end{example}

\begin{example}
    There are $\frac{1}{2}qr$ irreducible $G\times G$-modules $X_{[n]}\otimes X_{[m]}$ that contain $U_{a'}$, namely
    $$X_{[a]} \otimes X_{[sa'-a]}$$ for any irrep $X_{[a]}$.
    If $q$ is odd and $n'$ is even, there are $\frac{1}{4}(rs^2-6rs+10r)$ irreducible $G\times G$-modules $X_{[n]}\otimes X_{[m]}$ that contain $X_{[n']}$.\\
    If $q$ is odd and $n'$ is also odd, there are $\frac{1}{4}(rs^2-6rs+8r)$ irreducible $G\times G$-modules $X_{[n]}\otimes X_{[m]}$ that contain $X_{[n']}$. \\
    If $q$ is even, there are $\frac{1}{4}(rs^2-6rs+9r)$ irreducible $G\times G$-modules $X_{[n]}\otimes X_{[m]}$ that contain $X_{[n']}$.
\end{example}


\bibliographystyle{plain}
\bibliography{2025EDandMVP-bibliography}{}

\begin{thebibliography}{10}

\bibitem{MR2319497}
L.~Aburto, R.~Johnson, and J.~Pantoja.
\newblock The complex linear representations of {${\rm GL}(2,k)$}, {$k$} a
  finite field.
\newblock {\em Proyecciones}, 25(3):307--329, 2006.

\bibitem{MR4794599}
C.~Blanco~Villacorta, I.~Pacharoni, and J.~A. Tirao.
\newblock Matrix spherical functions on finite groups.
\newblock {\em J. Fourier Anal. Appl.}, 30(5):Paper No. 48, 37, 2024.

\bibitem{MR1450744}
R.~Camporesi.
\newblock The spherical transform for homogeneous vector bundles over
  {R}iemannian symmetric spaces.
\newblock {\em J. Lie Theory}, 7(1):29--60, 1997.

\bibitem{MR4162267}
T.~Ceccherini-Silberstein, F.~Scarabotti, and F.~Tolli.
\newblock {\em Gelfand triples and their {H}ecke algebras---harmonic analysis
  for multiplicity-free induced representations of finite groups}, volume 2267
  of {\em Lecture Notes in Mathematics}.
\newblock Springer, Cham, [2020] \copyright 2020.
\newblock With a foreword by Eiichi Bannai.

\bibitem{MR0954385}
R.~Gangolli and V.~S. Varadarajan.
\newblock {\em Harmonic analysis of spherical functions on real reductive
  groups}, volume 101 of {\em Ergebnisse der Mathematik und ihrer Grenzgebiete
  [Results in Mathematics and Related Areas]}.
\newblock Springer-Verlag, Berlin, 1988.

\bibitem{MR4748964}
J.C. Gardiner and S.P. Humphries.
\newblock Strong {G}elfand pairs of {${\rm SL}(2,p^n)$}.
\newblock {\em Comm. Algebra}, 52(8):3269--3281, 2024.

\bibitem{MR0052444}
R.~Godement.
\newblock A theory of spherical functions. {I}.
\newblock {\em Trans. Amer. Math. Soc.}, 73:496--556, 1952.

\bibitem{GuptaHassain}
A.~Gupta and M.~Hassain.
\newblock Tensor product of irreducible characters of
  $\mathrm{GL}_{2}(\mathbb{F}_q)$.
\newblock {\em Journal of Algebra and Its Applications}, 0(0):2650113, 0.

\bibitem{MR1313912}
G.~Heckman and H.~Schlichtkrull.
\newblock {\em Harmonic analysis and special functions on symmetric spaces},
  volume~16 of {\em Perspectives in Mathematics}.
\newblock Academic Press, Inc., San Diego, CA, 1994.

\bibitem{KnopRohrle}
F.~Knop and G.~R\"ohrle.
\newblock Spherical subgroups in simple algebraic groups.
\newblock {\em Compos. Math.}, 151(7):1288--1308, 2015.

\bibitem{MR4642865}
G.~Pezzini and M.~van Pruijssen.
\newblock On the extended weight monoid of a spherical homogeneous space and
  its applications to spherical functions.
\newblock {\em Represent. Theory}, 27:815--886, 2023.

\bibitem{P-S1983}
I.~Piatetski-Shapiro.
\newblock {\em Complex representations of {${\rm GL}(2,\,K)$}\ for finite
  fields {$K$}}, volume~16 of {\em Contemporary Mathematics}.
\newblock American Mathematical Society, Providence, RI, 1983.

\bibitem{Serre1977LinRep}
J.-P. Serre.
\newblock {\em Linear representations of finite groups}, volume Vol. 42 of {\em
  Graduate Texts in Mathematics}.
\newblock Springer-Verlag, New York-Heidelberg, french edition, 1977.

\bibitem{SimpsonFrameSutherland}
W.A. Simpson and J.~Sutherland Frame.
\newblock The character tables for {${\rm SL}(3,\,q)$}, {${\rm
  SU}(3,\,q\sp{2})$}, {${\rm PSL}(3,\,q)$}, {${\rm PSU}(3,\,q\sp{2})$}.
\newblock {\em Canadian J. Math.}, 25:486--494, 1973.

\bibitem{MR0774056}
D.~Stanton.
\newblock Orthogonal polynomials and {C}hevalley groups.
\newblock In {\em Special functions: group theoretical aspects and
  applications}, Math. Appl., pages 87--128. Reidel, Dordrecht, 1984.

\bibitem{Steinberg1951-Irreps}
R.~Steinberg.
\newblock The representations of {${\rm GL}(3,q), {\rm GL} (4,q), {\rm PGL}
  (3,q)$}, and {${\rm PGL} (4,q)$}.
\newblock {\em Canad. J. Math.}, 3:225--235, 1951.

\bibitem{Timashev}
D.A. Timashev.
\newblock {\em Homogeneous spaces and equivariant embeddings}, volume 138 of
  {\em Encyclopaedia of Mathematical Sciences}.
\newblock Springer, Heidelberg, 2011.
\newblock Invariant Theory and Algebraic Transformation Groups, 8.

\bibitem{MR0340399}
D.~Travis.
\newblock Spherical functions on finite groups.
\newblock {\em J. Algebra}, 29:65--76, 1974.

\bibitem{MR4815523}
P.~Turek.
\newblock Multiplicity-free induced characters of symmetric groups.
\newblock {\em Trans. Amer. Math. Soc.}, 377(12):8817--8876, 2024.

\end{thebibliography}

\end{document}